\documentclass[12pt]{amsart}

\usepackage{amscd}
\usepackage{amsfonts}
\usepackage{amssymb}
\usepackage{amsmath}
\usepackage{color}

\numberwithin{equation}{section}
\numberwithin{equation}{subsection}

\newtheorem*{pro*}{Proposition A}
\newtheorem*{thmA*}{Theorem A}
\newtheorem*{thmB*}{Theorem B}
\newtheorem*{thmD*}{Theorem D}

\newtheorem{ppro}{Proposition}[subsection]
\newtheorem{llemma}{Lemma}[subsection]
\newtheorem{ddefinition}{Definition }[subsection]
\newtheorem{ttheorem}{Theorem}[subsection]
\newtheorem{rremark}{Remark}[subsection]
\newtheorem{ccorollary}{Corollary}[subsection]

\newcommand{\al}{\alpha}
\newcommand{\be}{\beta}
\newcommand{\ga}{\gamma}

\newcommand{\se}{\theta}
\newcommand{\om}{\omega}
\newcommand{\we}{\wedge}
\newcommand{\lam}{\lambda}
\newcommand{\sig}{\sigma}
\newcommand{\ra}{{\rightarrow}}
\newcommand{\lra}{{\longrightarrow}}

\newcommand{\ie}{{i.e$.$\,}}
\newcommand{\cf}{{cf$.$\,}}

\newcommand{\CC}{{\mathbb C}}

\newcommand{\HH}{{\mathbb H}}
\newcommand{\KK}{{\mathbb K}}

\newcommand{\PP}{{\mathbb P}}
\newcommand{\RR}{{\mathbb R}}

\newcommand{\ZZ}{{\mathbb Z}}

\newcommand{\GL}{\mathop{\rm GL}\nolimits}

\newcommand{\Sim}{\mathop{\rm Sim}\nolimits}

\newcommand{\Conf}{\mathop{\rm Conf}\nolimits}
\newcommand{\Ker}{\mathop{\rm Ker}}

\newcommand{\dev}{\mathop{\rm dev}\nolimits}

\newcommand{\Diff}{\mathop{\rm Diff}\nolimits}

\begin{document}
\baselineskip 13pt 
\thispagestyle{empty}

\title[\ ]
{Conformally Lorentz parabolic structure and Fefferman Lorentz
metrics}

\author[\ ]{Yoshinobu KAMISHIMA}
\address{Department of Mathematics, Tokyo Metropolitan University,\linebreak
Minami-Ohsawa 1-1, Hachioji, Tokyo 192-0397, Japan.}
\email{kami@tmu.ac.jp}

\thanks{}
\date{\today}


\begin{abstract} We study conformal Fefferman-Lorentz manifolds of $(2n+2)$-dimension.
In order to do so, we introduce \emph{Lorentz parabolic structure} on
$(m+2)$-dimensional manifolds as a $G$-structure. By using causal conformal vector
fields preserving that structure, we shall establish two theorems on
compact Fefferman-Lorentz manifolds: One is the coincidence of vanishing curvature
between {\em Weyl conformal curvature tensor} of
Fefferman metrics on a Lorentz manifold $S^1\times N$ and {\em
Chern-Moser curvature tensor} on a strictly pseudoconvex
$CR$-manifold $N$. Another is the analogue of the conformal rigidity
theorem of Obata and Ferrand to the compact Fefferman-Lorentz
manifolds admitting noncompact closed causal conformal
transformations.
\end{abstract}

\maketitle

\setcounter{section}{0}
\section{Introduction}
A pseudo-Riemannian metric $g$ of signature $(m+1,1)$ is called a
Lorentz metric on an $m+2$-dimensional smooth manifold. An
$m+2$-dimensional Lorentz manifold $M$ is a smooth manifold equipped
with a Lorentz metric. Two Lorentz metrics $g,g'$ are conformal if
there exists a positive function $u$ on $M$ such that $g'=u\cdot g$.
The equivalence class $[g]$ of $g$ is a conformal class and
$(M,[g])$ is called a conformal Lorentz manifold. On the other hand,
given a $(2n+1)$-dimensional strictly pseudoconvex $CR$-manifold
$N$, C. Fefferman constructed a Lorentz metric on the product
$S^1\times N$ which has the following properties (\cf \cite{lee}):
\begin{itemize}
\item The conformal class of this metric is a $CR$-invariant.
\item $S^1$ acts as \emph{lightlike} isometries.
\end{itemize}
 It is interesting to know how the Fefferman-Lorentz
metric of $S^1\times N$ interacts on the $CR$-structure of $N$.
Compare \cite{lee} and general references
therein for the relation between Fefferman metrics and the Cartan connection.\\

In this paper, we shall take a different approach to the Fefferman-
Lorentz metrics  on  $(2n+2)$-dimensional manifolds by introducing a
$G_\CC$-structure called \emph{Fefferman-Lorentz parabolic}
structure. (Compare Section \ref{sec:G-structure}.)

Recall that \emph{conformal Lorentz structure} on an
$m+2$-dimensional manifold $M$ is an $\mathrm{O}(m+1,1)\times \RR^+$
-structure \cite{kobayashi}. An integrable $\mathrm{O}(m+1,1)\times
\RR^+$ -structure is \emph{conformally flat Lorentz structure} on
$M$. (See also \cite[p.10]{kobayashi}.) We focus on \emph{parabolic
subgroup} of $\mathrm{O}(m+1,1)$ which is defined as follows. Let
$\mathrm{PO}(m+1,1)$ be the real hyperbolic group, then the minimal
parabolic subgroup is isomorphic to either $\mathrm{O}(m+1)$ or the
similarity subgroup ${\Sim}(\RR^m)$, which lifts to the parabolic
subgroup of $\mathrm{O}(m+1,1)$ as $\mathrm{O}(m+1)\times \ZZ_2$ or
${\Sim}(\RR^m)\times \ZZ_2$ respectively.

For $m=2n$, we define a subgroup $G_\CC$ of ${\rm
Sim}(\RR^{2n})\times \RR^+\leq\mathrm{O}(2n+1,1)\times \RR^+$. If
${\rm U}(n+1,1)$ denotes the unitary Lorentz group which embeds into
${\rm O}(2n+2,2)$, then $G_\CC$ is characterized as the intersection
${\rm U}(n+1,1)\cap ({\rm Sim}(\RR^{2n})\times \RR^+)$ where ${\rm
Sim}(\RR^{2n})$ is identified with  the stabilizer ${\rm
PO}(2n+1,1)_\infty$ at the point at infinity.

Then we see that a Fefferman-Lorentz manifold admits a
Fefferman-Lorentz parabolic structure by reinterpreting  the proof
of \cite[(5.17) Theorem]{lee}.

One of our results concerns the relation between a $CR$-manifold $N$
and a Fefferman-Lorentz manifold $S^1\times N$. We provide
 a geometric proof to the coincidence of \emph{vanishing} between \emph{Weyl conformal
curvature tensor} and \emph{Chern-Moser curvature tensor}. This
result may be obtained as a special case of more general
calculations by Fefferman \cite{fe}.

\begin{thmA*}[Theorem \ref{equivalence-flatness}]\label{CR=Weyl}
A Fefferman-Lorentz manifold $S^1\times N$ is conformally flat
 if and only if $N$ is a spherical $CR$-manifold.
\end{thmA*}

By our definition, we obtain the following classes of conformally
flat Lorentz parabolic manifolds. We shall give compact examples in
each class. (See Section \ref{examples1}.)\begin{itemize}
\item Lorentz flat space forms.
\item Fefferman-Lorentz parabolic manifolds (locally modelled on\\ $(\hat{\rm U}(n+1,1),S^{2n+1,1})$).
\item Fefferman-Lorentz manifolds $S^1\times N$ where $N$ is a spherical $CR$-manifold.
\end{itemize}

In the second part, we study the \emph{Vague Conjecture}
 \cite{dam-gro} that the existence of a global geometric flow  determines a compact geometric manifold
uniquely, \ie \emph{isomorphic to the standard model
with flat $G$-structure}. The celebrated theorem of Obata and Ferrand provides
a supporting example for this,
\ie if a closed group $\RR$ acts
conformally on a compact Riemannian manifold, then it is conformal to the standard sphere $S^n$.

We study the analogue of the theorem of Obata and Ferrand to
compact Lorentz manifolds. In general, it is not true only by the existence
of a noncompact conformal closed subgroup $\RR$.
It is a problem which conformal group gives rise to
an affirmative answer to the compact Lorentz case,
\ie a compact Lorentz manifold is conformal to the Lorentz model $S^{n-1,1}$.
C. Frances and K. Melnick gave a
sufficient condition on the nilpotent dimension of
a nilpotent Lie group acting conformally on a compact Lorentz manifold.
(Compare \cite{cfk} more generally for pseudo-Riemannian manifolds.)

We prove affirmatively the theorem of Obata and Ferrand to
 the compact Fefferman-Lorentz manifolds under the existence of
 \emph{two dimensional causal abelian Lie groups}.
Let $M=S^1\times N$ be a compact Fefferman-Lorentz
manifold on which $S^1$ acts as Lorentz lightlike isometries. Denote
by $\mathcal C_{{\Conf}(M,g)}(S^1)$ the centralizer of $S^1$ in
${\Conf}(M,g)$.

\begin{thmB*}[Theorem \ref{th:o-l}]\label{th:vague}
Suppose that
  $\mathcal C_{{\Conf}(M,g)}(S^1)$ contains a closed
noncompact subgroup of dimension $1$ at least. Then $M$ is
conformally equivalent to the two-fold cover $S^1\times S^{2n+1}$ of
the standard Lorentz manifold $S^{2n+1,1}$.
\end{thmB*}
\emph{Dimension $2$} for the dimension of $\mathcal
C_{{\Conf}(M,g)}(S^1)$ is rather small relative to the nilpotent
dimension in \cite{cfk}. However the fact that $S^1$ of $M$ is the
group of lightlike isometries with respect to our Fefferman-Lorentz
metric is more geometric. In fact, let $N_{-}$ be a Lorentz
hyperbolic $3$-manifold $\mathrm {PSL}(2,\RR)/\Gamma$. Then
$M=S^1\times N_{-}$ admits a conformally flat Lorentz metric on
which $S^1$ acts as \emph{spacelike} isometries. (See Remark
\ref{rem:1-theorem8.12}). The metric is not a Fefferman-Lorentz
metric but there exists a \emph{two dimensional noncompact abelian
Lie group} $S^1\times \RR$ acting isometrically on $M$. On the other
hand, if we note that $N_{-}$ admits a spherical $CR$-structure,
then $M$ does admit a Fefferman-Lorentz metric whose conformal
Lorentz group is compact (Lorentz isometry group). In addition, the
above group $S^1\times \RR$ is not a conformal group of the
underlying Fefferman-Lorentz metric.

\tableofcontents

\section{Preliminaries}
\subsection{Cayley-Klein model} We start with the Cayley-Klein projective model for
$\KK=\RR,\CC$ or $\HH$. (Compare
\cite{bcdgm},\cite{ck},\cite{cf},\cite{kam1},\cite{on},\cite{wo} for
instance.) There is an equivariant principal bundle:
\[
\KK^*\ra ({\rm GL}(n+2,\KK)\cdot \KK^*,\KK^{n+2}-\{0\}) \stackrel
P{\lra} ({\rm PGL}(n+2,\KK), \KK\PP^{n+1}).
\] Fix nonnegative integers $p,q$ such that $n=p+q$.
The nondegenerate Hermitian form on $\KK^{n+2}$ is defined by
\begin{equation}\label{q-inner} \mathcal
B(x,y) = {\bar x}\sb 1y\sb 1+\cdots +
 {\bar x}\sb{p+1} y\sb{p+1}-{\bar x}\sb{p+2} y\sb{p+2}-\cdots-
 {\bar x}\sb{n+2} y\sb{n+2}
\end{equation}
Denote by ${\rm O}(p+1,q+1;\KK)$ the subgroup of ${\GL}(n+2, \KK)$:
\begin{equation*}
\begin{split}
\{A\in \mathop{\GL}(n+2, \KK)\ |\ \mathcal B(Ax,Ay)=\mathcal B(x,y),
x,y\in \KK^{n+2}\}.
\end{split}
\end{equation*} The group ${\rm O}(p+1,q+1;\KK)$ leaves invariant the
$\KK$-cone in $\KK^{n+2}-\{0\}$:
\begin{equation}\label{k-cone}
\begin{split}
V \sb 0 &=\{x \in \KK^{n+2}-\{0\}\ \vert \ \mathcal B(x,x)=0\}.
\end{split}
\end{equation}
Put
\begin{equation}\begin{split}\label{3-models}
P(V\sb 0)=\left\{\begin{array}{lcr}
S^{p,q}& \approx &S^{q}\times S^{p}/\ZZ_2\\
S^{2p+1,2q} & \approx &S^{2q+1}\times S^{2p+1}/S^1\\
S^{4p+3,4q}&\approx &S^{4q+3}\times S^{4p+3}/S^3
\end{array}\right.
\end{split}
\end{equation}
Let ${\rm PO}(p+1,q+1;\KK)$ denote the image of ${\rm
O}(p+1,q+1;\KK)$ in ${\rm PGL}(n+2,\KK)$. According to $\KK=\RR,\CC$
or $\HH$, we have the \emph{nondegenerate} flat geometry of
signature $(p,q)$.
\[{\small
\left\{\begin{array}{lr}
({\rm PO}(p+1,q+1),S^{p,q}) &\ \ \mbox{Conformally flat geometry\ \ \ \  \ \ \ \ \ \ \ \ \ \ \ }\\
({\rm PU}(p+1,q+1),S^{2p+1,2q}) &\ \ \mbox{Spherical $CR$ geometry\
\ \ \ \ \ \
\ \ \ \ \  \ \ \ \ \ \ \ }\\
({\rm PSp}(p+1,q+1),S^{4p+3,4q}) &\mbox{Flat pseudo-conformal q$CR$
geometry}
\end{array}\right.
}\] (Compare \cite{kam1}.) 
In particular, when $p=n,q=0$, \ie positive definite case, we have
the usual horospherical geometry, \ie the geometry on the boundary
of the real, complex or quaternionic hyperbolic spaces. When $p=n-1,
q=1$,  $({\rm PO}(n,2),S^{n-1,1})$ is said to be the
\emph{$n$-dimensional conformally flat Lorentz geometry}.

\subsection{Causality}  Let $\xi$ be a vector field on a Lorentz manifold
$(M,g)$. We recall \emph{causality of vector fields} (\cf
\cite{on}).

\begin{ddefinition}[Causal vector fields]
\label{def:causalfields} Let $x\in M$.
\begin{equation}
\left\{\begin{array}{lcr}
\xi\ \mbox{is spacelike} & g(\xi_x,\xi_x)>0 & \mbox{whenever}\ \xi_x\neq 0.\\
\xi\ \mbox{is lightlike} & g(\xi_x,\xi_x)=0 & \mbox{whenever}\ \xi_x\neq 0.\\
\xi\ \mbox{is timelike}  & g(\xi_x,\xi_x)<0 & \mbox{whenever}\ \xi_x\neq 0.\\
\end{array}\right.
\end{equation}
\end{ddefinition}

Each vector $\xi_x$ is called a \emph{causal vector}. Suppose that
$\xi$ is a vector field defined on a domain $\Omega$ of $M$. If
$\xi_x\neq 0$ and $\xi_x$ is a causal vector at each point $x\in
 \Omega$, then $\xi$ is said to
be a causal vector field on $\Omega$.

\section{Conformally Lorentz parabolic geometry}\label{sec:Logeometry}

\subsection{Conformally  Lorentz parabolic structure}\label{sec:G-structure}
 Let $\{e_1,\dots,e_{m}\}$ be the
standard orthonormal basis of $\RR^{m+2}$ with respect to the
Lorentz inner product $\mathcal B$; $\mathcal
B(e_i,e_j)=\delta_{ij}$ $(1\leq i,j\leq m+1)$, $\mathcal
B(e_{m+2},e_{m+2})=-1$. (See \eqref{q-inner} for $\KK=\RR$.) Set
\[
\ell_1={e_1+e_{m+2}}{/}\sqrt 2, \ \ell_{m+2}={e_1-e_{m+2}}{/}\sqrt
2.\] Putting $\mathcal B=\langle \ , \ \rangle$ on $\RR^{m+2}$,
$\mathsf{F}=\{\ell_1, e_2,\dots,e_{m+1},\ell_{m+2}\}$ is a new basis
such that $\langle
\ell_1,\ell_1\rangle=\langle\ell_{m+2},\ell_{m+2}\rangle=0$,\,
$\langle\ell_1,\ell_{m+2}\rangle= 1$. The symmetric matrix
$\mathsf{I}_{m+1}^{1}$ with respect to this basis $\mathsf{F}$
is described as
\begin{equation}\label{uni}
\mathsf {I}_{m+1}^{1}=\left(\begin{array}{c|c|c}
0 &  0     \cdots  0       & 1\\
\hline
0 &     & 0 \\
\vdots & \mathrm{I}_{m} &\vdots\\
 0 &     & 0\\
 \hline
1 &  0 \cdots 0 & 0 \\
\end{array}\right).
\end{equation}
Note that
\[
\mathrm{O}(m+1,1)=\{A\in\mathrm{GL}(m+2,\RR)\,|\ A\,
\mathsf{I}_{m+1}^{1}{}^tA=\mathsf{I}_{m+1}^{1}\}.
\]

The similarity subgroup is described in $\mathrm{O}(m+1,1)$ as
follows:

 {\small
\begin{equation}\label{similaritygroup}
{\Sim}(\RR^m)=\left\{\left(\begin{array}{ccc}
\lambda & x           &-\mbox{\Large $\frac{|x|^2}{2\lambda}$}\\
  & B           &-\mbox{\Large $\frac{B\,{}^tx}{\lambda}$}\\
 \mbox{{\Large $0$}} &                     & \mbox{\Large$\frac{1}{\lambda}$} \\
\end{array}\right)\, \mbox{\huge$|$} \begin{array}{cc} \lam\in \RR^+, & B\in\mathrm{O}(m)\\
                                                        x\in\RR^{m}  &
                                                         \\
                                                         \end{array}
\right\}.
\end{equation}
}
We introduce the following subgroup in ${\GL}(m+2,\RR)$:
\begin{equation}\label{similo}
 G_\RR={\Sim}(\RR^m)\times \RR^+.
\end{equation}
Let $P=\sqrt u\cdot Q\in G_\RR$ $(Q\in{\Sim}(\RR^m), \sqrt
u\in\RR^+)$. For the Lorentz inner product: $$\langle
x,y\rangle=x\,\mathsf{I}_{m+1}^{1}{}^ty,$$
$(x=(x_1,\dots,x_{m+2}),y=(y_1,\dots,y_{m+2})\in \RR^{m+2}$) as
above,
\begin{equation*}
\langle xP,yP\rangle= u\cdot xQ\,\mathsf{I}_{m+1}^{1}{}^tQ{}^ty =
u\cdot x\,\mathsf{I}_{m+1}^{1}{}^ty=u\cdot\langle x,y\rangle,
\end{equation*}
\ie a $G_\RR$-structure defines a conformal class of Lorentz metrics
on an $(m+2)$-manifold.

 Let $G_\CC\leq {\GL}(2n+2,\RR)$ be a subgroup defined by
{\small
\begin{equation}\label{similoF}
G_\CC=\left\{\left(\begin{array}{ccc}
u & \sqrt u x           &-\mbox{\Large $\frac{|x|^2}2$}\\
  & \sqrt u B           &-B\,{}^tx\\
\mbox{{\large $0$}}  &                     & 1 \\
\end{array}\right)\, \mbox{\huge$|$}  u\in \RR^+, x\in\CC^{n}, B\in\mathrm{U}(n)
\right\}.
\end{equation}
} If $P$ is an element of $G_\CC$ as above, then
\begin{equation}\label{eq:parabolic}
\begin{split}
P&=\sqrt u\cdot Q,\\
Q=&\left(\begin{array}{ccc}
\sqrt u &  x           &-\mbox{\Large $\frac{|x|^2}{2\sqrt{u}}$}\\
  & B           &-\mbox{\Large$\frac1{\sqrt u}$}B\,{}^tx\\
\mbox{{\Large $0$}}  &                     & \mbox{\Large$\frac1{\sqrt u}$} \\
\end{array}\right).
\end{split}\end{equation}
It is easy to see that $Q$ is an element of
$\mathrm{Sim}(\RR^{2n})\leq\mathrm{O}(2n+1,1)$. So $G_\CC$ is a
subgroup of $G_\RR$ for $m=2n$. A homomorphism $P\mapsto Q$ gives
rise to an isomorphism of $G_\CC$ onto
$\mathrm{Sim}(\CC^n)=\CC^{n}\rtimes(\mathrm{U}(n)\rtimes \RR^+)$ of
$\mathrm{Sim}(\RR^{2n})$.

\begin{ddefinition}\label{FL}
\par\ \par
\begin{itemize}
\item A $G_\RR$-structure on an $(m+2)$-manifold is called
\emph{conformally  Lorentz parabolic} structure.
 An $(m+2)$-manifold
is said to be a conformally Lorentz parabolic manifold if it admits
a $G_\RR$-structure.

 \item A $G_\CC$-structure
 on a $(2n+2)$-manifold is \emph{Fefferman-Lorentz parabolic structure}. In other words, the
Fefferman-Lorentz parabolic structure is a
 reduction of $G_\RR$ to $G_\CC$.
A Fefferman-Lorentz parabolic manifold is a $(2n+2)$-dimensional
manifold equipped with a Fefferman-Lorentz parabolic structure.
\end{itemize}
\end{ddefinition}

\begin{ddefinition}\label{fl-structure}
Let $M$ be a Fefferman-Lorentz parabolic manifold. Then
$\mathop{\Conf}_{\rm FLP}(M)$ is the group of conformal
transformations preserving the Fefferman-Lorentz parabolic
structure.
\end{ddefinition}

\subsection{Lorentz similarity geometry $(\mathcal L{\rm sim}(\RR^{m+2}),\RR^{m+2})$}
\label{L-similarity}
 Let
$\RR^{m+2}$ be the $(m+2)$-dimensional euclidean space  equipped
with a Lorentz inner product (\cf Section \ref{sec:G-structure}).
Form the Lorentz similarity subgroup $\mathcal L{\rm
sim}(\RR^{m+2})=\RR^{m+2}\rtimes ({\rm O}(m+1,1)\times \RR^+)$ from
the affine group ${\rm Aff}(\RR^{m+2})=\RR^{m+2}\rtimes {\rm
GL}(m+2,\RR)$. If an $(m+2)$-manifold $M$ is locally modelled on
$\RR^{m+2}$ with coordinate changes lying in $\mathcal L{\rm
sim}(\RR^{m+2})$, then $M$ is said to be a \emph{Lorentz similarity}
manifold.

\begin{ppro}\label{Lsimilaritymanifold}
Let $M$ be an $(m+2)$-dimensional compact Lorentz similarity
manifold with virtually solvable fundamental group. Then $M$ is
either a Lorentz flat parabolic manifold or finitely covered by a
Hopf manifold $S^{m+1}\times S^1$, or an $m+2$-torus $T^{m+2}$.
\end{ppro}

\begin{proof}
Given a compact Lorentz similarity manifold $M$, there exists a
developing pair $\displaystyle (\rho,\dev):(\pi_1(M),\tilde M)\ra
(\mathcal L{\rm sim}(\RR^{m+2}),\RR^{m+2})$. Suppose that $\pi_1(M)$
is virtually solvable. 
Let $L:\mathcal L{\rm sim}(\RR^{m+2})\ra {\rm O}(m+1,1)\times \RR^+$
be the linear holonomy homomorphism. Then a subgroup of finite index
in $L(\rho(\pi_1(M))$ is solvable in ${\rm O}(m+1,1)\times \RR^+$ so
it belongs to a maximal amenable subgroup which is either isomorphic
to $({\rm O}(m+1)\times {\rm O}(1))\times \RR^+$ or to
$({\Sim}(\RR^m)\times \ZZ_2)\times \RR^+$ up to conjugate. By
Definition \ref{FL}, the latter case implies that $M$ (or two
fold-cover) is a Lorentz flat parabolic manifold.
 If a subgroup of finite index in
$L(\rho(\pi_1(M))$ lies in $({\rm O}(m+1)\times {\rm O}(1))\times
\RR^+$, then $M$ is a similarity manifold where ${\rm O}(m+1)\times
{\rm O}(1)\leq {\rm O}(m+2)$. It follows from the result by Fried
that $M$ is covered finitely by an $m+2$-torus $T^{m+2}$ or a Hopf
manifold $S^{m+1}\times S^1$.

\end{proof}

 The Lorentz similarity geometry
contains \emph{Lorentz flat} geometry\\ $({\rm E}(m+1,1),\RR^{m+2})$
where ${\rm E}(m+1,1)=\RR^{m+2}\rtimes {\rm O}(m+1,1)$. If $M$ is an
$m+2$-dimensional compact Lorentz flat manifold, then it is known
that $M$ is geodesically complete and the fundamental group of a
compact complete Lorentz flat manifold is virtually solvable.
Applying the above proposition, we have
\begin{ccorollary}
Any $m+2$-dimensional compact Lorentz flat manifold is finitely
covered by an $m+2$-torus or an infrasolvmanifold.
\end{ccorollary}
\begin{proof}The holonomy homomorphism $\rho:
\pi_1(M)\ra \RR^{m+2}\rtimes {\rm O}(m+1,1)$ reduces to a
homomorphism: $\rho:\pi_1(M)\ra \RR^{m+2}\rtimes {\Sim}(\RR^m)\times
\ZZ_2$. Put $\Gamma=\rho(\pi_1(M))$ which is a virtually solvable
discrete subgroup. As ${\Sim}(\RR^m)=\RR^m\rtimes ({\rm O}(m)\times
\RR^+)$, a subgroup $\Gamma'$ of finite index in $\Gamma$ is
conjugate to a discrete subgroup of a solvable Lie group
$G=\RR^{m+2}\rtimes (\RR^m\rtimes (T^m\times \RR^+))$. Let
$M'=\rho^{-1}(\Gamma')\backslash\tilde M$ be a finite covering of
$M$. As $\mathop{\dev}:\tilde M\ra \RR^{m+2}$ is a diffeomorphism,
it follows that $M'\cong \Gamma'\backslash \RR^{m+2}=
\Gamma'\backslash G/H$ where $H=\RR^m\rtimes (T^m\times \RR^+)$.

\end{proof}

\subsection{Fefferman-Lorentz manifold}\label{fl-def}
 In \cite{f} Fefferman has shown that when
$N$ is a $(2n+1)$-dimensional \emph{strictly pseudoconvex}
$CR$-manifold, $S^1\times  N$ admits a Lorentz metric $g$
on which $S^1$ acts as lightlike isometries.
 We recall the construction of the  metric from \cite{lee}.
Let $({\rm Ker}\,\om, J)$ be a $CR$-structure on $N$
with characteristic (Reeb) vector field $\xi$ for some contact form
$\om$. The circle $S^1$ generates the vector field $\mathcal S$ on
$S^1\times N$ (extending trivially on $N$). Note that
\begin{equation}\label{eq:decomp}
T(S^1\times N)=\langle\mathcal S\rangle \oplus\langle\xi\rangle
\oplus {\rm Ker}\,\om.
\end{equation}
Let $\displaystyle ({\rm Ker}\,\om)\otimes \CC=\{Y_1,\dots,
Y_{n}\}\oplus \{\bar Y_{1},\dots, \bar Y_{n}\}$ be the canonical
decomposition for $J$ for which we choose such as $d\om(Y_i,\bar
Y_{j})=-\mathbf {i}\delta_{ij}$. As usual, letting $\displaystyle
X_i=Y_i+\bar Y_i/\sqrt 2$, $\displaystyle
X_{n+i}=\mathbf{i}(Y_i-\bar Y_i)/\sqrt 2$, it implies that
$\displaystyle ({\rm Ker}\,\om)=\{ X_1,\dots, X_{2n}\}$  such that
$JX_i=X_{n+i}$ and
\begin{equation}\label{normalframe}
d\om(JX_i,X_{j})=\delta_{ij}\ \ (i=1,\dots,n).
\end{equation}
 This gives a (real) frame $\{\mathcal S,\xi, X_1,\dots,X_{2n}\}$ at a
neighborhood of $S^1\times N$.
 Let $\se^i$ be the dual frame to $X_i$ $(i=1,\dots,2n)$. From
\eqref{normalframe}, note that
\begin{equation}\label{eq:dualX}
d\om(J-,-)=\sum_{i=1}^{2n}\se^i\cdot\se^i\ \, \mbox{on}\ {\rm
Ker}\,\om.
\end{equation}

Let $\displaystyle {dt}$ be a $1$-form on $S^1\times N$ such that
\begin{equation}\label{dtvalue}
{dt}(\mathcal S)=1, \, \  {dt}(V)=0 \ \ (\forall\, V\in TN).
\end{equation} Put $\displaystyle\eta=\om^1\we\cdots \we\om^n$ where
$\displaystyle\om^\al=\se^\al+\mathbf{i}\se^{n+\al}$. By Proposition
(3.4) of \cite{lee} there exists a unique real $1$-form $\sigma$ on
$S^1\times N$ satisfying that
\begin{equation}\label{eq:dualT}
\begin{split}
 d(\om\wedge \eta)&=\mathbf{i}(n+2)\sig\we\om\we\eta,\\
\sigma\we d\eta\we\bar\eta&=\mathrm {Tr}(d\sigma)\mathbf
{i}\sigma\we \om\we\eta\we\bar\eta.\\
\end{split}
\end{equation}
 The explicit form of $\sigma$ is obtained from \cite[(5.1) Theorem]{lee} that
\begin{equation}\label{leescomp}
\sigma=\frac 1{n+2}\left(dt+\mathbf{i}P^*\om_\al^{\al}-\frac
1{2(n+1)}\rho\cdot P^*\om\right).
\end{equation}Here $P:S^1\times N\to N$ is the canonical projection and
$\om_\al^{\be}$ is a connection form of $\om$
such that
\begin{equation*}
\begin{split}
d\om&=\mathbf{i}\delta_{\al\be}\om^\al\we\om^{\bar\be},\\
d\om^\al&=\om^\be\we\om_\be^\al+\om\we\tau^\be.
\end{split}
\end{equation*}
The function $\rho$ is the Webster scalar curvature on $N$. (Since
we chose $h_{\al\bar\be}=\delta_{\al\be}$, note that $\displaystyle
-\frac{\mathbf{i}} 2h^{\al\bar\be}dh_{\al\bar\be}=0$ in the equation
(5.3) of \cite[ (5.1) Theorem]{lee}.) It follows that
\begin{equation}\label{sigvavalue}
\sigma(\mathcal S)=\frac 1{n+2}.
\end{equation}

 Define a symmetric $2$-form
$$\sigma\odot\om=\sigma\cdot\om+\om\cdot\sigma.$$ Extending $\se^i(\mathcal S)=0$
and $\om(\mathcal S)=0$, we have a Fefferman - Lorentz metric on
$S^1\times N$ by
\begin{equation}\label{lometric}
\begin{split}
g(X,Y)&=\sigma(X)\cdot \om(Y)+\om(X)\cdot\sigma(Y)+d\om(JX^{h},Y^{h})\\
 &=\sigma\odot \om(X,Y)+\sum_{i=1}^{2n}\se^i\cdot\se^i(X,Y).
\end{split}
\end{equation}
Here $X^h$ stands for the horizontal part of $X$, \ie $X^h\in
\mathrm{Ker}\,\om$. By \eqref{sigvavalue}
we have that
\begin{equation}\label{pairofLo}
g(\xi,\mathcal S)=\frac 1{n+2}.
\end{equation}
Since $g(\mathcal S,\mathcal S)=0$, $g$ becomes a Lorentz metric on
$S^1\times N$. In particular $S^1$ acts as lightlike isometries of
$g$.

 The following result has been achieved by
Lee \cite{lee}.
 We shall give an elementary proof of \emph{invariance}
 from the viewpoint of $G$-structure.

\begin{ttheorem}\label{th:CRequivL}
A strictly pseudoconvex $CR$-structure on $N$ gives a conformal
class of Lorentz metrics on $S^1\times N$. Indeed $S^1\times N$
admits a Fefferman-Lorentz parabolic structure $(G_\CC$-structure$)$.
\end{ttheorem}

\begin{proof}
 Suppose that $({\rm Ker}\,\om', J)$
represents the same $CR$-structure on $N$. Then it follows that
$\om'=u\cdot \om$ for some positive function $u$ on $N$. Let
$\{\mathcal S',\xi', X'_1,\dots,X'_{2n}\}$ be another frame on the
neighborhood of $S^1\times N$ for $\om'$. Since $\mathcal S'$
generates the same $S^1$, note that
\begin{equation}\label{eq:V}
\mathcal S=\mathcal S'.
 \end{equation}There exist $\ x_i\in\RR$ $(i=1,\dots,2n)$ for which the
characteristic vector field $\xi'$ is described as
\begin{equation}\label{eq:chravec}
\xi=u\cdot \xi'+x_1\sqrt uX'_1+\cdots+x_{2n}\sqrt uX'_{2n}.
\end{equation}
 As $u\cdot d\om=d\om'$ and $d\om(J-,J-)=d\om(-,-)$ on
${\rm Ker}\,\om$, there exists a $B=({b_{i}}^{j})\in\mathrm{U}(n)$
such that
\begin{equation}\label{eq:unitary}
X_i=\sqrt u \sum_{k}^{}{b_{i}}^{k}X'_k.
\end{equation}
Two frames $\{\mathcal S,\xi, X_1,\dots,X_{2n}\}$, $\{\mathcal
S',\xi', X'_1,\dots,X'_{2n}\}$ are uniquely determined each other by
the equations \eqref{eq:V}, \eqref{eq:chravec}, \eqref{eq:unitary}.

It suffices to prove that another Lorentz metric $g'$ is conformal
to $g$:
$$ g'=
\sigma'\odot\om'+\sum_{i=1}^{2n}{\se'}^i\cdot{\se'}^i.$$ The
equations \eqref{eq:V}, \eqref{eq:chravec}, \eqref{eq:unitary}
determine the relation between the dual frames $\{\om,
\se^1,\dots,\se^{2n},\sigma\}$, $\{\om',
{\se'}^1,\dots,{\se'}^{2n},\sigma'\}$.
\begin{equation}\label{dualrel}
\begin{split}
\om'&=u\cdot\om,\\
{\se'}^i&=\sqrt u\sum_{j} {b_{j}}^{i}\se^{j}+\sqrt ux_i\cdot \om,\\
\end{split}\end{equation}
Moreover, by the uniqueness property of $\sigma$ from
\eqref{eq:dualT}, $\sigma'$ is transformed into the following form
(\cf \cite[(5.16) Proposition]{lee}):
\begin{equation}\label{dualre2}
\sigma'= \sigma-\sum_{i,j}^{}{b_{j}}^{i}x_i\se^j-\frac{|x|^2}2\om.
\end{equation}

Using \eqref{eq:parabolic}, the above equations show that {\small
\begin{equation}\begin{split}
 (\om',{\se'}^1,\dots,
{\se'}^{2n},\sigma')&=(\om,{\se}^1,\dots, {\se}^{2n},\sigma)
\left(\begin{array}{ccc}
u & \sqrt u x           &-\mbox{\Large $\frac{|x|^2}2$}\\
  & \sqrt u B           &-B\,{}^tx\\
 \mbox{{\large $0$}}  &                     & 1 \\
\end{array}\right)\\
&=(\om,{\se}^1,\dots, {\se}^{2n},\sigma)\sqrt u\cdot Q.
\end{split}\end{equation}
} As $B\in \mathrm{U}(n)$ and $x\in \CC^n$ because the basis
$\{X_1,\dots, X_{2n}\}$ invariant under $J$, note that $\sqrt u\cdot
Q\in G_\CC$ (\ie  $Q\in \mathrm{Sim}_{\CC}(n)$).  Hence the
$CR$-structure defines a Fefferman-Lorentz parabolic structure on
$M$ (\cf Definition \ref{FL}). Moreover, a calculation shows

\begin{equation}\label{lorentz2}
\begin{split}
g'&=\sigma'\odot\om'+d\om'(J- ,- )=
\sigma'\cdot\om'+\om'\cdot\sigma'+\sum{\se'}^i\cdot{\se'}^i\\
&=(\om',{\se'}^1,\dots, {\se'}^{2n},\sigma')\, \mathsf
{I}_{2n+1}^{1}{}^t(\om',{\se'}^1,\dots, {\se'}^{2n},\sigma')\\
&=(\sqrt u)^2(\om,{\se}^1,\dots, {\se}^{2n},\sigma)\,Q
\mathsf{I}_{2n+1}^{1} {}^tQ\,
{}^t(\om,{\se}^1,\dots, {\se}^{2n},\sigma)\\
&=u\cdot (\om,{\se}^1,\dots, {\se}^{2n},\sigma)\,\mathsf
{I}_{2n+1}^{1}\,{}^t(\om,{\se}^1,\dots, {\se}^{2n},\sigma)\\
&=u(\sigma\odot\om+\sum{\se}^i\cdot{\se}^i)=u\cdot g,\\
\end{split}\end{equation}
Hence the $CR$-structure determines a conformal class of
Fefferman-Lorentz metric $g$.

\end{proof}

\section{Conformally flat Fefferman-Lorentz manifold}\label{classi1}
\subsection{Confomally flat Fefferman-Lorentz model}\label{flat-nil}

 Let \[
 V \sb
0=\{x=(x_1,\dots,x_{2n+4})\in \RR^{2n+4}-\{0\}\ |\ \mathcal
B(x,x)=0\}\] be as in \eqref{k-cone} for $\KK=\RR$. In this case,
when the Hermitian bilinear form is defined by
\begin{equation}\label{c-form}
 \langle z,w\rangle =\bar z_1w_1+\cdots+\bar z_{n+1}w_{n+1}-\bar z_{n+2}
w_{n+2}\ \ \mbox{on}\ \CC^{n+2},
\end{equation} $V_0$ is identified with
\[
\{z=(z_1,\dots,z_{n+2})\in \CC^{n+2}-\{0\}\ |\ \langle z,z\rangle
=0\}.\]  Let ${\rm U}(n+1,1)$ be the unitary Lorentz group with the
center $S^1$. Obviously the two-fold cover of $S^{2n+1,1}$ is
contained in $V_0$, \ie $S^1\times S^{2n+1}\subset V_0$ but not
invariant under $\mathrm{U}(n+1,1)$. Consider the commutative
diagram. {\footnotesize
\begin{equation}\label{good daigram}
\begin{CD}
(\ZZ_2,\RR^*) @= (\ZZ_2, \RR^*) \\
@VVV       @VVV   \\
(S^1,\CC^*) @>>> \CC^{n+2}-\{0\} @>P_\CC>> \CC\PP^{n+1}\\
{||}@.   \bigcup@. \bigcup\\
(S^1,\CC^*) @>>> S^1\times S^{2n+1}\subset V_0 @>P_\CC>> S^{2n+1}=P_\CC(V_0)\\
@VVV  @V{P_\RR}VV    {||} \\
S^1/\ZZ_2=S^1@>>>  S^{2n+1,1}=P_\RR(V_0) @>P>> S^{2n+1}=P(S^{2n+1,1})\\
{||}@.  \bigcap@.          \bigcap \\
S^1 @>>> \RR\PP^{2n+3} @>{P}>> \CC\PP^{n+1}.
\end{CD}
\end{equation}
}
 Put
$$\hat{\rm U}(n+1,1)= {\rm U}(n+1,1)/\ZZ_2$$
where $\ZZ_2$ is a cyclic group of order two in $S^1$. The natural
embedding ${\rm U}(n+1,1)\ra {\rm O}(2n+2,2)$ induces an embedding
of Lie groups:
 \[ \hat{\rm U}(n+1,1)\ra {\rm PO}(2n+2,2).
 \]

Let \begin{equation}\label{sphericalform}\begin{split}
\om_0=-\mathbf{i}\sum_{j=1}^{n+1}\bar z_jdz_j\, ,\   \ \
\xi=\mathop{{\sum}}_{j=1}^{n+1}(x_j\frac d{dy_j}\ -\ y_j\frac
d{dx_j})
\end{split}
\end{equation} be the standard contact form on $S^{2n+1}$ with
the characteristic vector field $\xi$. Note that
$\displaystyle\omega_0(\xi)=\mathop{{\sum}}_{j=1}^{n+1}|z_j|^2=1$ on
$S^{2n+1}$. As in Section \ref{fl-sp}, $\om_0$ is the connection
form on the principal bundle : $\displaystyle S^1\ra
S^{2n+1}\stackrel{\pi}\lra \CC\PP^n$. The spherical $CR$-structure
$(\mathrm{Ker}\, \om_0,J_0)$ defines a Lorentz metric on $S^1\times
S^{2n+1}$:

\begin{equation}\label{f-lorentz}
g^0(X,Y)=\sigma_0\odot P^*_\CC\om_0(X,Y)+
d\om_0(J_0{P_\CC}_*X,{P_\CC}_*Y)
\end{equation}
where $\sigma_0$ is obtained from \eqref{leescomp} (\cf
\eqref{sigmaform}).

\begin{ppro}\label{confflatlo}
The group ${\mathrm{U}}(n+1,1)$ acts conformally on $S^1\times
S^{2n+1}$ with respect to $g^0$. Especially, so does
$\hat{\mathrm{U}}(n+1,1)$ on $(S^{2n+1,1},\hat g^0)$.
\end{ppro}

\begin{proof}
 Let
$\mathrm{U}(n+1,1)=(\mathrm{U}(n+1)\times \mathrm{U}(1))\cdot
(\mathcal N\times \RR^+)$ be the Iwasawa decomposition in which
there is the equivariant projection:
\begin{equation}\label{p-map}
\begin{CD}
({\rm U}(n+1,1),V_0)@>(P,P_\CC)>>({\rm PU}(n+1,1),S^{2n+1}).
\end{CD}
\end{equation}
If $t_\se=e^{\mathbf{i}\theta}\in \mathrm{ZU}(n+1,1)$ which is the
center $S^1$ of $\mathrm{U}(n+1,1)$, then  by the form
\eqref{leescomp}  it follows that
\begin{equation}\label{sigmainv} t_\se^*dt=dt, \ \ \, t_\se^*\sigma_0=\sigma_0.
\end{equation} Let $\mathrm{U}(n+1)$
be the maximal compact subgroup of $\mathrm{PU}(n+1,1)$. If
$\ga=P(\tilde \ga)\in \mathrm{U}(n+1)$, then $\ga^*\om_0=\om_0$ from
\eqref{sphericalform}. Then $\tilde \ga^* g^0=g^0$, \ie $\tilde \ga$
acts as an isometry of $S^1\times
S^{2n+1}$.\\

Suppose that $\ga=P(\tilde \ga)\in \mathcal N\times \RR^+$ where
$\mathcal N$ is the Heisenberg Lie group such that $\mathcal
N\cup\{\infty\}=S^{2n+1}$. Recall from Section $2$ of \cite{kam2}
that if $(t,(z_1,\cdots,z_n))$ is the coordinate of $\mathcal
N=\RR\times\CC^n$, then the contact form $\omega_{\mathcal N}$ on
$\mathcal N$ is described as:
\begin{equation}\label{nil-con}
\omega_{\mathcal N}
=dt+\mathop{{\sum}}_{j=1}^{n}(x_jdy_j-y_jdx_j)=dt+{\rm Im} \langle
z,dz\rangle.
\end{equation}(Here $\displaystyle \langle
z,w\rangle =\mathop{\sum}_{i=1}^{n} \bar z_iw_i$ and ${\rm Im}\, x$
is the imaginary part of $x$.) An element
 $g=((a,z),\lambda\cdot A)\in\mathcal N\rtimes
 (\mathrm{U}(n)\times \RR^+)$
acts on $(t,w)\in \mathcal N$ as
\begin{equation}\label{cr-action}
g\cdot (t,w)= (a+\lambda^2\cdot t-{\rm Im}\langle z,\lambda\cdot
A\cdot w\rangle, z+\lambda\cdot A\cdot w).
\end{equation}In particular,
when $\ga=(a,z)\in\mathcal N$, then $\ga^*\om_\mathcal
N=\om_\mathcal N$, while $\RR^+=\langle\ga_\theta\rangle$ satisfies
that $\displaystyle \gamma_\theta(t,z)=(e^{2\theta}\cdot
t,e^{\theta}\cdot z)$
 on $\mathcal N$ so $${\gamma_\theta}^*\omega_\mathcal
N=d(e^{2\theta}\cdot t)+ {\rm Im}\langle e^{\theta}\cdot z,
d(e^{\theta}\cdot z)\rangle=e^{2\theta}\cdot\omega_\mathcal N.$$

As $(\om_0,J)$ and $(\om_\mathcal N,J)$ define the same spherical
$CR$-structure on $\mathcal N$, there exists a smooth function $u$
on $\mathcal N$ such that $\om_{\mathcal N}=u\cdot \om_0$. Let
\[
g_\mathcal N=\sigma_\mathcal N\odot P^*\om_\mathcal N+ d\om_\mathcal
N(J_0P_*-,P_*-)\] be the Lorentz metric on $S^1\times \mathcal N$
where $P:S^1\times \mathcal N\ra \mathcal N$ is the projection. Then
it follows from Theorem \ref{th:CRequivL} that

\begin{equation}\label{flatg0}
g_\mathcal N=u\cdot g^0.
\end{equation}
As above, it is easy to check that if $P(\tilde\ga)=\ga\in\mathcal
N$, then $\ga^*g_\mathcal N=g_\mathcal N$ and if
$P(\gamma_\theta)=\gamma_\theta\in \RR^+$, then $\ga_\se^*g_\mathcal
N=e^{2\theta}\cdot g_\mathcal N$. (Note that the equation
${\gamma_\theta}^*\omega_\mathcal N=e^{2\theta}\cdot\omega_\mathcal
N$ implies that ${\gamma_\theta}^*d\omega_\mathcal
N=e^{2\theta}\cdot d\omega_\mathcal N$ because $e^{2\theta}$ is
constant.) As a consequence, if $\ga=P(\tilde \ga)\in \mathcal
N\times \RR^+$, then there exists a positive constant $\tau$ such
that $\tilde \ga^*g_\mathcal N=\tau\cdot g_\mathcal N$. Letting a
positive function $v=\ga^*u^{-1}\cdot\tau\cdot u$ on $\mathcal N$,
it is easy to see that $$\tilde\ga^*g^0=v\cdot g^0\ \ \mbox{on}\
\,S^1\times \mathcal N.$$
 Since this is true on a neighborhood at any point in
$S^1\times S^{2n+1}$, $\tilde \ga$ acts conformally on $S^1\times
S^{2n+1}$ with respect to $g^0$. For  either $\tilde
\ga\in\mathrm{U}(n+1)\times \mathrm{U}(1)$ or $\tilde\ga\in \mathcal
N\times \RR^+$, the above observation shows that every element of
$\mathrm{U}(n+1,1)$ acts as conformal transformation on $S^1\times
S^{2n+1}$ with respect to the Lorentz metric $g^0$.

\end{proof}

\medskip
Note from \eqref{flatg0} that  $g_\mathcal N=u\cdot g^0$.
The Weyl conformal curvature tensor satisfies that
$W(g^0)=W(g_\mathcal N)$ on $S^1\times \mathcal N$. In order to
prove that the Fefferman-Lorentz metric $g^0$ is a conformally flat
metric, we calculate the Weyl conformal curvature tensor of
$g_\mathcal N$ on $S^1\times \mathcal N$ directly.

In view of the contact form $\om_\mathcal N$ on $\mathcal N$ (\cf
\eqref{nil-con}), $d\om_\mathcal N(J_0-,-)$ is the euclidean metric
$\displaystyle \hat g_\CC=2\mathop{{\sum}}_{j=1}^{n}|dz_j|^2$ on
$\CC^n$. As $\om_\mathcal N$ is the connection form of the principal
bundle: $\mathcal R\to \mathcal N\stackrel{\pi}\lra\CC^n$, it
follows from \eqref{leescomp} that
$$
\sigma_\mathcal N=\frac
1{n+2}\left(dt+\mathrm{i}P^*\pi^*\varphi_\al^{\al}- \frac
1{2(n+1)}\pi^*\rho\cdot P^*\om_\mathcal N\right).$$ Since
$(\CC^n,\hat g_\CC)$ is flat, it follows that $\varphi_\al^{\al}=0$,
 $\rho=\pi^*s=0$, which shows
\begin{equation}\label{sigma=dt}
\sigma_\mathcal N=\frac 1{n+2}dt.
\end{equation}
The metric $g_\mathcal N$ reduces to the following:
\[g_\mathcal N=\frac 1{n+2}dt\odot
P^*\om_\mathcal N+P^* \hat g_\CC.\]

It is easy to check that the following are equivalent:
\begin{itemize}
\item[{\bf (i)}] $X\in\mathfrak
C^{\perp}$ where $\mathfrak C=\langle \mathcal S,\xi\rangle$ induced
by $S^1\times \mathcal R$.
\item[{\bf (ii)}] $g_\mathcal N(X,\mathcal S)=0$ and $g_\mathcal
N(X,\xi)=0$.
\item[{\bf (iii)}] $X\in P^*\mathop{\Ker}\, \om_\mathcal N$,\, $dt(X)=0$.
\end{itemize} As a consequence,
$\mathfrak C^{\perp}=\mathop{\Ker}\, \om_\mathcal N$. Putting
$\mathcal C=S^1\times \mathcal R$,
there is a pseudo-Riemannian submersion:
\begin{equation}\label{subnil}
\begin{CD} \mathcal C @>>> (S^1\times \mathcal N,g_\mathcal N)
@>{\pi}>>(\CC^n,\hat g_\CC).
\end{CD}
\end{equation}
The Riemannian curvature tensor $\hat R$ on the flat space $\CC^n$
is zero,
$$\hat g_\CC(\hat R(\pi_*X,\pi_*Y)\pi_*Z,\pi_*W)=0.$$

If $X,Y\in P^*\mathop{\Ker}\, \om_\mathcal N$, then $[X,Y]^{\mathcal
V}\in\langle \xi\rangle$ where $\xi$ is the characteristic vector
field for $\om_\mathcal N$. (Here $X^{\mathcal V}$ stands for the
fiber component of the vector $X$.) Since
$dt(\xi)=(n+2)\sigma_\mathcal N(\xi)=0$ from \eqref{sigma=dt},
$$g_\mathcal N(\xi,\xi)=0,\, \ \mbox{\ie}\ \xi\, \ \mbox{is lightlike}.$$

We apply the O'Neill's formula (\cf \cite[(3.30)]{c-e} for example)
to the pseudo-Riemannian submersion of \eqref{subnil}:
\begin{equation}\label{oneill}
\begin{split}
&g_\mathcal N(R(X,Y)Z,W)=\hat g_\CC(\hat R(\pi_*X,\pi_*Y)\pi_*Z,\pi_*W)\\
&\ \ \ \ \ \   \ \ \ \ \ +\frac14g_\mathcal N([X,Z]^{\mathcal
V},[Y,W]^{\mathcal V})
-\frac14 g_\mathcal N([Y,Z]^{\mathcal V},[X,W]^{\mathcal V})\\
&\ \ \ \ \ \ \ \ \ \ \ \ \ \ +\frac12 g_\mathcal N([Z,W]^{\mathcal
V},[X,Y]^{\mathcal V}]). \end{split}
\end{equation}
This shows that
\begin{llemma}\label{curvaturecal1}
\begin{equation*}
R_{XYZW}=g_\mathcal N(R(X,Y)Z,W)=0\, \ (\forall\, X,Y,Z,W\in
P^*\mathop{\Ker}\, \om_\mathcal N).
\end{equation*}
\end{llemma}
\begin{llemma}\label{curvaturecal2}
\begin{equation*}
R_{\xi ABC}=0\ \ (\forall\, A,B,C\in T(S^1\times \mathcal N)).
\end{equation*}

\end{llemma}

\begin{proof}
Put $\om=\om_\mathcal N$, $g=g_\mathcal N$ and
$\mathop{\Ker}\, \om=P^*\mathop{\Ker}\, \om_\mathcal N$.
Let $\nabla$ be a
covariant derivative for $g$ on $S^1\times \mathcal N$;
\begin{equation*}\begin{split}
2g(\nabla_XY,Z)&=Xg(Y,Z)+Yg(X,Z)-Zg(X,Y)\\
 & +g([X,Y],Z)+g([Z,X],Y)+g([Z,Y],X).\\
 \end{split}
 \end{equation*}As ${\rm Ker}\, \om$ is $\mathcal C$-invariant,
 we note the following.
 \begin{equation}\label{iv}
 [X,\xi]=[X,\mathcal S]=0\ \ (\forall\, X\in {\rm Ker}\, \om).
\end{equation}
Put $\displaystyle \sigma=\sigma_\mathcal N=\frac 1{n+2}dt$. By {\bf
(iii)},
\begin{equation}\label{v}
\sigma([X,Y])=0\ \ (\forall\, X,Y\in {\rm Ker}\, \om).
\end{equation} We may choose $X,Y,Z$ to be orthonormal vector fields
in ${\rm Ker}\, \om$. Then
\begin{equation*}\begin{split}
2g(\nabla_X\xi,Z)&=g([Z,X],\xi)=\sigma([Z,X])=0 \ \
(\forall\,X,Z\in{\rm Ker}\, \om).
 \end{split}
 \end{equation*}
 Similarly from \eqref{iv},
\begin{equation*}\begin{split}
2g(\nabla_X\xi,\mathcal S)&=Xg(\xi,\mathcal S)=X(\frac 1{n+2})=0,\\
2g(\nabla_X\xi,\xi)&=0.\\
 \end{split}
 \end{equation*}This implies that
\begin{equation}\label{vi}
\nabla_X\xi=0.
\end{equation}It follows similarly that
\begin{equation}\label{geodesicSxi}
2g(\nabla_\mathcal S\xi,Z)=0,\ \, 2g(\nabla_\mathcal S\xi,\mathcal
S)=0,\ \, 2g(\nabla_\mathcal S\xi,\xi)=0.
\end{equation}This shows that

\begin{equation}\label{vii}
\nabla_\mathcal S\xi=0.
\end{equation}
It is easy to see that $\nabla_\mathcal S\mathcal S=\nabla_\xi
\xi=0$, \ie the orbits of $S^1$ and $\RR$ are geodesics. From these,
we obtain that
\begin{equation}\label{viii}
\nabla_\xi A=0\, \ (\forall\, A\in T(S^1\times \mathcal N)).
\end{equation}
This implies that
 $R_{\xi ABC}=0$ $(\forall\, A,B,C\in T(S^1\times \mathcal N))$.

\end{proof}

We set formally
\begin{equation}\label{compatible}
J\mathcal S=0,\ J\xi=0
\end{equation}so that $J$ is defined on
$T(S^1\times \mathcal N)=\{\mathcal S,\xi\}\oplus {\rm Ker}\, \om$.

\begin{llemma}\label{curvaturecal3}
\begin{equation*}
\nabla_\mathcal SA=-\frac 1{n+2}JA.
\ \ (\forall\, A\in T(S^1\times \mathcal N)).
\end{equation*}
\end{llemma}

\begin{proof}
For the vector field $\mathcal S$, we see that
\begin{equation}\label{ix}
2g(\nabla_X\mathcal S,\xi)=0,\ \ 2g(\nabla_X\mathcal S,\mathcal
S)=0.
\end{equation}
As we assumed that $g(X,Z)=0$ and $g(X,X)=g(Z,Z)=1$, calculate
\begin{equation*}\begin{split}
2g(\nabla_X\mathcal S,Z)&=\om([Z,X])\cdot\sigma(S)
=\frac 1{n+2}\cdot\om([Z,X])\\
&=-\frac 2{n+2}\cdot d\om(Z,X)=
-\frac 2{n+2}\cdot \hat g_\CC(Z,\hat JX)\\
&=-\frac 2{n+2}\cdot g(JX,Z).\\
\end{split}
\end{equation*}Using \eqref {ix}, we obtain that
\begin{equation}\label{perppart}
\nabla_X\mathcal S=-\frac 1{n+2}JX.
\end{equation}
As $[\mathcal S,X]=0$, note that
$$\nabla_\mathcal SX=-\frac 1{n+2}JX.$$
Since $\nabla_\mathcal S\xi=\nabla_\mathcal S S=0$ from \eqref
{vii}, we have that $\displaystyle \nabla_\mathcal SA=-\frac
1{n+2}JA$ \ $(\forall\, A\in T(S^1\times \mathcal N))$.

\end{proof}

\begin{llemma}\label{curvaturecal4}
The remaining curvature tensor $R_{ABCD}$ on $S^1\times \mathcal N$
becomes as follows.
\begin{enumerate}
\item[(1)] $\displaystyle R_{\mathcal SX\mathcal S Y}=
-\frac 1{(n+2)^2}\cdot g_\mathcal N(X,Y)$\ \ $(\forall\, X,Y\in
P^*\mathop{\Ker}\, \om_\mathcal N)$.
\item[(2)] $R_{\mathcal SXYZ}=0$\ \ $(\forall\, X,Y,Z\in
P^*\mathop{\Ker}\, \om_\mathcal N)$.
\end{enumerate}
\end{llemma}

\begin{proof}

Using Lemma \ref{curvaturecal3},
\begin{equation*}\begin{split}
R(X,\mathcal S)\mathcal
S&=\nabla_X\nabla_\mathcal S(\mathcal S)-\nabla_\mathcal
S\nabla_X(\mathcal
S)-\nabla_{[\mathcal S,X]}\mathcal S\\
&=-\nabla_\mathcal S(-\frac 1{n+2}JX) =\frac 1{(n+2)^2}X.
\end{split}\end{equation*}

It follows that

\begin{equation*}\label{onepoint}
\begin{split}
R_{\mathcal SX\mathcal SY}&=-R_{X\mathcal S\mathcal SY}\\
&=-g(R(X,\mathcal S)\mathcal
S,Y)=-\frac 1{(n+2)^2}g(X,Y)
\end{split}\end{equation*}
As $X,Y$ are orthonormal and $\sigma([X,Y])=0$ by \eqref{v}, it
follows that
\begin{equation*}
2g(\nabla_XY,\xi)=g([X,Y],\xi)=0,
\end{equation*}so there exists a function $a(X,Y)$
such that
\begin{equation*}
\nabla_XY\equiv a(X,Y)\xi \ \, {\rm mod}\,{\rm Ker}\, \om.
\end{equation*}In particular it follows that
\begin{equation}\label{mod}
\nabla_XJY-J\nabla_XY\equiv a(X,JY)\xi  \ \, {\rm mod}\,{\rm Ker}\, \om.
\end{equation}

Let $\hat\nabla$ be a covariant derivative for the K\"ahler metric
$\hat g_\CC$ on $\CC^n$ as before. Recall from \cite[(3.23)
p.67]{c-e} that
\begin{equation}\label{perp-metric}
g(\nabla_XY,Z)=\hat g_\CC(\hat\nabla_{\pi_*X}\pi_*(Y),\pi_*Z)\
\ (\forall\, X,Y,Z\in
\mathop{\Ker}\, \om).
\end{equation}If we note that
the complex structure $\hat J$ is parallel with respect to
$\hat g_\CC$, \ie $\hat\nabla\hat J=\hat J\hat\nabla$, then
\begin{equation*}
\begin{split}
g(\nabla_XJY,Z)&=\hat g_\CC(\hat\nabla_{\pi_*X}\pi_*(JY),\pi_*Z)\\
&=\hat g_\CC(\hat\nabla_{\pi_*X}\hat J\pi_*(Y),\pi_*Z)\\
&=\hat g_\CC(\hat J\hat\nabla_{\pi_*X}\pi_*(Y),\pi_*Z)\\
&=-\hat g_\CC(\hat\nabla_{\pi_*X}\pi_*(Y),\pi_*JZ)\\
&=-g(\nabla_XY,JZ)=g(J\nabla_XY,Z).\\
\end{split}\end{equation*}
\eqref {mod} implies that
\begin{equation}\label{exactly}
\nabla_XJY-J\nabla_XY=a(X,JY)\xi.
\end{equation}
As $[X,\mathcal S]=0$,  Lemma \ref{curvaturecal3} shows that
\begin{equation*}
\begin{split}
R(X,\mathcal S)Y&=\nabla_X\nabla_\mathcal S Y-
\nabla_\mathcal S\nabla_X Y\\
&=-\frac 1{n+2}(\nabla_XJY-J\nabla_XY)\\
&=-\frac 1{n+2}\, a(X,JY)\xi.\\
\end{split}
\end{equation*}By (1) of Lemma \ref{curvaturecal4}, it follows that
\begin{equation*}
\begin{split}
g(R(X,\mathcal S)Y,\mathcal S)&=-\frac 1{(n+2)^2}a(X,JY)\\
&=R_{X\mathcal SY\mathcal S}=-\frac 1{(n+2)^2}g(X,Y).
\end{split}
\end{equation*}Hence
$a(X,JY)=g(X,Y)$ so that we obtain
\begin{equation}\label{sxyz}
R(X,\mathcal S)Y=-\frac 1{n+2}g(X,Y)\xi.
\end{equation}
It follows that
\begin{equation*}\label{rsxyz}
R_{\mathcal S XYZ}=-R_{X\mathcal S YZ}= -g(R(X,\mathcal
S)Y,Z)=0.\end{equation*}

\end{proof}

\begin{ppro}\label{Weyl1}Let $(S^1\times \mathcal N,g_\mathcal
N)$ be a Fefferman-Lorentz nilmanifold of dimension $2n+2$. Then the
following hold.

\begin{itemize}
\item[(1)] The scalar curvature function $S=0$.
\smallskip
\item[(2)] The Ricci tensor has the following form.
\end{itemize}

\begin{itemize}
\item [{\rm (i)}] $R_{YZ}=0$ \ \ $(\forall\,Y,Z\in {\rm
Ker}\,\om)$.
\item[{\rm (ii)}] $R_{\xi A}=0$ \ $(\forall\,A\in T(S^1\times \mathcal N))$.
\item[{\rm (iii)}] $R_{\mathcal SY}=0$  \ $(\forall\,Y\in {\rm Ker}\,\om)$.
\item[{\rm (iv)}]  $R_{\mathcal S\mathcal S}=-2n/(n+2)^2$.
\end{itemize}

\end{ppro}

\begin{proof}

Note that $g(\xi,\xi)=g(\mathcal S,\mathcal S)=0$. Then
\begin{equation*}\begin{split}
S&=R_{ABCD}g^{AC}g^{BD}\\
&=R_{XYZW}g^{XZ}g^{YW}+ R_{\xi BCD}g^{\xi C}g^{BD}+R_{\mathcal SBCD}g^{\mathcal S C}g^{BD}\\
\end{split}
\end{equation*}where $X,Y,Z,W\in {\rm Ker}\, \om$ and $B,C,D\in T(S^1\times \mathcal N)$.
The first term is zero; $R_{XYZW}=0$ by Lemma \ref{curvaturecal1}.
The second term $R_{\xi BCD}=0$ by Lemma \ref{curvaturecal2}.
According to whether $g^{\mathcal S C}=0$ or $g^{\mathcal S
\xi}={n+2}$, the third term becomes $R_{\mathcal SBCD}g^{\mathcal S
C}g^{BD}=R_{\mathcal SB\xi D}g^{\mathcal S \xi}g^{BD}=0$ because
$R_{\mathcal SB\xi D}=0$ by Lemma \ref{curvaturecal2} again. Hence
the scalar curvature $S=0$.

The Ricci tensor satisfies that
\begin{equation*}
\begin{split}
R_{\mathcal SY}&=R_{A\mathcal SBY}g^{AB}\\
&=R_{X\mathcal SZY}g^{XZ}+R_{\xi \mathcal SBY}g^{\xi B}
+R_{A\mathcal S\xi Y}g^{A\xi}+ R_{A\mathcal S\mathcal
SY}g^{A\mathcal S}.
\end{split}
\end{equation*}Then $R_{X\mathcal SZY}=-R_{\mathcal SXZY}=0$ by (2)
of Lemma \ref{curvaturecal4}.
$R_{\xi\mathcal SBY}=0, R_{A\mathcal S\xi Y}=R_{\xi Y A\mathcal
S}=0$ by Lemma \ref{curvaturecal2}.  When $g^{\xi\mathcal S }=n+2$,
$R_{\xi\mathcal S\mathcal SY}=0$ as above. According to whether
$g^{A\mathcal S}=0$ or $g^{\xi\mathcal S }=n+2$, the third term
$R_{A\mathcal S\mathcal SY}g^{A\mathcal S}=0$. So $R_{\mathcal
SY}=0$, (iii) follows.  Similarly (i), (ii) follow by calculations:
\begin{equation*}
\begin{split}
R_{YZ}&=R_{AYBZ}g^{AB}=R_{\xi YBZ}g^{\xi B}+
 R_{\mathcal S YBZ}g^{\mathcal SB}\\
&=R_{\mathcal SY\xi Z}g^{\mathcal S\xi}=0.
\end{split}\end{equation*}
\begin{equation*}
R_{\xi A}=R_{B\xi C A}g^{BC}=0 \ \, (\forall\,A\in T(S^1\times
\mathcal N)).
\end{equation*}
As we chose $g(X,X)=1$,  (1) of Lemma \ref{curvaturecal4} implies
that
\begin{equation*}
\begin{split}
R_{\mathcal S\mathcal S}&=R_{A\mathcal SB\mathcal S}g^{AB}\\
&=R_{X\mathcal SX\mathcal S}g^{XX}+ R_{X\mathcal S\xi\mathcal
S}g^{X\xi}+ R_{\xi\mathcal SB\mathcal
S}g^{\xi B}+R_{\mathcal S\mathcal SB\mathcal S}g^{\mathcal SB}\\
&=2n\cdot R_{X\mathcal SX\mathcal S}=-2n\cdot \frac 1{(n+2)^2}.
\end{split}\end{equation*}
 This shows (iv).

\end{proof}

Let $W_{ABCD}$ be the Weyl conformal curvature tensor of $g_\mathcal
N$ on $S^1\times \mathcal N$ of dimension $2n+2$. Recall that
\begin{equation*}\begin{split}
W_{ABCD}&=R_{ABCD}+\frac
1{2n}\left(R_{BC}g_{AD}-R_{BD}g_{AC}-R_{AC}g_{BD}+R_{AD}g_{BC}\right)\\
&\ \ \ \ +\frac S{2n(2n+1)}\left(g_{BD}g_{AC}-g_{BC}g_{AD}\right).
\end{split}\end{equation*}

\begin{ppro}\label{Lorentznilflat}
The Fefferman-Lorentz manifold $(S^1\times \mathcal N,g_\mathcal N)$
is a conformally flat Lorentz nilmanifold.
\end{ppro}

\begin{proof}
We shall prove that all the Weyl conformal curvature tensors vanish.
By (1) of Proposition \ref {Weyl1}, this reduces to
\begin{equation*}\begin{split}
W_{ABCD}&=R_{ABCD}+\frac
1{2n}\left(R_{BC}g_{AD}-R_{BD}g_{AC}-R_{AC}g_{BD}+R_{AD}g_{BC}\right).
\end{split}\end{equation*}

It follows from (i) of Proposition \ref{Weyl1} and Lemma
\ref{curvaturecal1},
\begin{equation}\label{perppart1}
W_{XYZW}=0 \ \ (\forall\,X,Y,Z,W\in {\rm Ker}\,\om).\\
\end{equation}
It follows from (i), (ii), (iii) of Proposition \ref{Weyl1} and
Lemma \ref{curvaturecal2},
\begin{equation}\label{perppart2}
\begin{split}
W_{\xi YZW}&=0 \ \ (\forall\,Y,Z,W\in {\rm Ker}\,\om),\\
W_{\xi Y\xi W}&=0,\ W_{\xi Y\xi \mathcal S}=0,\ W_{\xi Y\mathcal S
W}=0.
\end{split}
\end{equation}
Similarly by (ii), (iii) of Proposition \ref{Weyl1},
\begin{equation}\label{perppart3}
\begin{split}
W_{\xi \mathcal SZW}&=0,\\
W_{\xi \mathcal S\xi W}&=-W_{\xi W \xi\mathcal S}=0.\\
\end{split}
\end{equation}

Let \begin{equation*}\begin{split} W_{\xi\mathcal S\xi\mathcal S}&=
\frac 1{2n}\left(R_{\mathcal S\xi}g_{\xi\mathcal S }-R_{\mathcal
S\mathcal S}g_{\xi\xi}-R_{\xi\xi}g_{\mathcal S\mathcal
S}+R_{\xi\mathcal S }g_{\mathcal S\xi}\right).
\end{split}\end{equation*}Since  $R_{\mathcal
S\mathcal S}\neq 0$ but $g_{\xi\xi}=0$, it follows that
\begin{equation}\label{perppart4}
W_{\xi\mathcal S\xi\mathcal S}=0.
\end{equation}
As $\displaystyle W_{\xi\mathcal S\mathcal SW}= \frac
1{2n}\left(R_{\mathcal S\mathcal S}g_{\xi W}-R_{\mathcal S
W}g_{\xi\mathcal S}-R_{\xi\mathcal S}g_{\mathcal S W}+R_{\xi W}
g_{\mathcal S\mathcal S}\right)$ but $g_{\xi W}=0$, it follows that

\begin{equation}\label{perppart5}
W_{\xi\mathcal S\mathcal SW}=0, \ W_{\xi\mathcal S\mathcal S\xi}=0.
\end{equation}

From \eqref{perppart2}, \eqref{perppart3}, \eqref{perppart4},
\eqref{perppart5}, the Weyl tensors containing $\xi$ are zero. It
follows similarly
\begin{equation}\label{perppart6}
\begin{split}
W_{\mathcal S\xi Z W}&=-W_{\xi\mathcal S Z W}=0, \\
W_{\mathcal S\xi \xi W}&=-W_{\xi\mathcal S\xi W}=0,\\
W_{\mathcal S\xi\mathcal SW}&=-W_{\xi\mathcal S\mathcal SW}=0,\\
W_{\mathcal S\xi\xi\mathcal S}&=0, \ W_{\mathcal S\xi\mathcal
 S\xi}=0,\\
W_{\mathcal SY\xi W}&=W_{\xi W\mathcal SY}=0,\\
W_{\mathcal SY\xi\mathcal S}&=W_{\xi\mathcal S\mathcal SY}=0.
\end{split}
\end{equation}

By (1) of  Lemma \ref{curvaturecal4} and (iv) of Proposition
\ref{Weyl1}, we calculate
\begin{equation}\begin{split}
W_{\mathcal SY\mathcal SW}&=R_{\mathcal SY\mathcal SW}\\
&+\frac 1{2n}\left(R_{Y\mathcal S}g_{\mathcal SW}-R_{YW} g_{\mathcal
S\mathcal S}-R_{\mathcal S\mathcal S}g_{YW}+R_{\mathcal S W
}g_{Y\mathcal S}\right)\\
&=R_{\mathcal SY\mathcal SW}-\frac 1{2n}R_{\mathcal S\mathcal
S}g_{YW}\\
&=\frac {-1}{(n+2)^2}g(Y,W)-\frac 1{2n}\left(\frac
{-2n}{(n+2)^2}g(Y,W)\right)=0.
\end{split}\end{equation}
So all the terms containing $\mathcal S$ are zero. We thus conclude
that the Weyl conformal curvature tensors of $(S^1\times \mathcal
N,g_\mathcal N)$ vanish.

\end{proof}

\begin{ppro}
The Fefferman-Lorentz metric $g^0$ is a conformally flat Lorentz
metric on $S^1\times S^{2n+1}$.
\end{ppro}

\begin{proof}
If we note  that $g_\mathcal N=u\cdot g^0$ as before,
the Weyl conformal curvature tensor satisfies that
$W(g^0)=W(g_\mathcal N)=0$ on $S^1\times \mathcal N$. Since
${\mathrm{U}}(n+1,1)$ acts conformally and transitively on
$S^1\times S^{2n+1}$, $W(g^0)=0$ on $S^{2n+1,1}$.
\end{proof}

\subsection{Examples of Fefferman-Lorentz manifolds I}\label{fl-sp}
We shall give examples of conformally flat Fefferman-Lorentz
manifolds which admit causal Killing fields.
 Consider principal $S^1$-bundles as a connection bundle over a K\"ahler manifold $W$.
\[
S^1\ra N^{2n+1}\stackrel{\pi}\lra W.
\]
There exists a $1$-form $\om$ such that $d\om=\pi^*\Omega$ for which
$\Omega=\mathbf{i}\delta_{\al\be}\om^\al\we \om^{\bar \be}$ is the
K\"ahler form on $W$. As $\mathop{\Ker}\,\om$ is isomorphic to $TW$
at each point of $W$.
Let $J$ be a complex structure on
$\mathop{\Ker}\,\om$ obtained from that of $W$ by the pullback of
$\pi_*$. Then $(\mathop{\Ker}\,\om, J)$ is a strictly pseudoconvex
$CR$-structure on $N$. Let $\xi$ be a characteristic vector field
induced by $S^1$. Let $P\colon S^1\times N\to N$ be the projection
as before. We have a Lorentz metric on $S^1\times N$:
\begin{equation}\label{f-lorentz}
g=\sigma\odot P^*\om+d\om(JP_*-,P_*-).
\end{equation}
Let $d\om^\al=\om^\be\we \varphi_\be^\al$ be the structure equation
on $W$ for the K\"ahler form $\Omega$. As
$d\om=\mathbf{i}\delta_{\al\be}\pi^*\om^\al\we \pi^*\om^{\bar \be}$,
the structure equation for $\om$ becomes
\[
d\pi^*\om^\al=\pi^*\om^\be\we \pi^*\varphi_\be^\al. 
\]
Let $s$ be the scalar curvature of $W$. Then the Webster scalar
function is defined as $\rho=\pi^*s$. By the definition,
\begin{equation}\label{sigmaform}
\sigma=\frac
1{n+2}\left(dt+\mathrm{i}P^*\pi^*\varphi_\al^{\al}-\frac
1{2(n+1)}\pi^*s\cdot P^*\om\right).
\end{equation}

As $\xi$ is characteristic, $\om(\xi)=1$, and $P_*\xi=\xi$,
$\pi_*(\xi)=0$ on $W$.
 From \eqref{sigmaform}, we obtain that
\begin{equation}\label{sigmasimple}
\sigma(\xi)=-\frac 1{2(n+1)(n+2)}s.
\end{equation}
Let $c$ be a positive constant. When $W$ is the complex projective
space $\CC\PP^n$, a complex torus $T^n_\CC$ or a complex hyperbolic
manifold $\HH^n_\CC/\Gamma$ of constant holomorphic sectional
curvature $c$, $0$, $-c$ respectively, the scalar curvature
$\displaystyle s= \frac {n(n+1)c}2,\,0,\,-\frac {n(n+1)c}2$
respectively. It follows from \eqref{sigmasimple} that
\[\sigma(\xi)=-\frac {nc}{4(n+2)},\, 0
,\, \frac {nc}{4(n+2)}\]respectively.  We obtain that
\begin{equation}\label{g-simple}
g(\xi,\xi)=\left\{\begin{array}{lr}
 -\displaystyle{\frac {nc}{2(n+2)}}&\ \ \mbox{for}\ \CC\PP^n,\\
 & \\
 \ \ 0 & \mbox{for}\ T^n_\CC,\\
  & \\
 \displaystyle{ \frac {nc}{2(n+2)}}& \ \ \ \mbox{for}\
 \HH^n_\CC/\Gamma.
\end{array}\right.
\end{equation} Note that there are principal $S^1$-bundles as a
connection bundle:
\begin{equation}\label{3-connection}
\begin{CD}
S^1@>>> S^{2n+1}@>\pi>> \CC\PP^n,\\
\RR@>>> \mathcal N@>\pi>> \CC^n,\\
S^1@>>> V^{2n,1}_{-1}@>\pi>> \HH^n_\CC.\\
\end{CD}
\end{equation}Here $\mathcal N$ is the Heisenberg Lie group with
isometry group ${\rm Isom}(\mathcal N)=\mathcal N\rtimes {\rm
U}(n)$. There is a discrete subgroup $\Delta\leq\mathcal N\rtimes
{\rm U}(n)$ whose quotient has a principal fibration: $S^1\ra
\mathcal N/\Delta\lra T^n_\CC$. Using the complex coordinates,
$V^{2n,1}_{-1}$ is
defined as 
\begin{equation}\label{cover-hyeprbolic}
\{(z_1,\dots, z_{n+1})\in \CC^{n+1}\ |\
|z_1|^2+\cdots+|z_{n}|^2-|z_{n+1}|^2=-1\}.
\end{equation}
 Moreover, the complement $S^{2n+1}-S^{2n-1}$ is identified with
$V_{-1}^{2n,1}$. Then the group ${\rm U}(n,1)$ acts transitively on
$V_{-1}^{2n,1}$ whose stabilizer at a point is isomorphic to ${\rm
U}(n)$.
 Moreover, there exists a discrete cocompact subgroup $\Gamma$ of ${\rm U}(n,1)$
 such that the image $\pi(\Gamma)$ is a torsionfree discrete cocompact subgroup
 of ${\rm U}(n,1)={\rm Isom}(\HH^n_\CC)$.

Note that if the first summand $S^1$ of $S^1\times N$ generates a
vector field $\mathcal S$, then $\displaystyle \sigma(\mathcal
S)=\frac 1{n+2}dt(\mathcal S)=1$ and $P^*\om(\mathcal S)=
\om(P_*\mathcal S)=0$. Then $g(\mathcal S,\mathcal S)=0$, \ie $S^1$
is lightlike.

\begin{ppro}\label{3typeskilling}
Let $X$ be one of the compact conformally flat Fefferman-Lorentz
manifolds $S^1\times S^{2n+1}$, $S^1\times \mathcal N/\Delta$ or
$S^1\times V_{-1}^{2n,1}/\Gamma$. Then $X$ admits a lightlike
Killing vector field $\mathcal S$ and a timelike $($resp. lightlike,
spacelike$)$ Killing vector field $\xi$.
\end{ppro}

\section{Uniformization}\label{uni}
\subsection{Uniformization of Fefferman-Lorentz parabolic
 manifolds}
 We prove the following uniformization concerning
Fefferman-Lorentz parabolic manifolds.

\begin{ppro}\label{pro:flat-almost complex}
Let $M$ be a $(2n+2)$-dimensional Fefferman-Lorentz parabolic
 manifold $($\ie admits a $G_\CC$-structure.$)$ If $M$ is conformally
flat, then $M$ is uniformizable with respect to $(\hat {\rm
U}(n+1,1),S^{2n+1,1})$.
\end{ppro}

\begin{proof} Suppose that $M$ is a
conformally flat Lorentz $(2n+2)$-manifold. By the definition, there
exists a collection of charts
$\{U_\al,\varphi_\al\}_{\al\in\Lambda}$. Let $\varphi_\al: U_\al\ra
S^{2n+1,1}$, $\varphi_\be:U_\be\ra S^{2n+1,1}$ be charts with
$U_\al\cap U_\be\neq\emptyset$. The coordinate change
$\varphi_\al\circ \varphi_\be^{-1}:\varphi_\be(U_\al\cap U_\be)\ra
\varphi_\al(U_\al\cap U_\be)$ extends to a transformation
$g_{\al\be}\in {\rm PO}(2n+2,2)$ of $S^{2n+1,1}$.

By the existence of $G_\CC$-structure, we have the principal frame
bundle: $\displaystyle G_\CC\ra P\lra M$ where $G_\CC\leq{\rm
O}(2n+1,1)\times \RR^+$. This bundle restricted to each neighborhood
$U_\al$ gives the trivial principal bundle on each neighborhood of
$S^{2n+1,1}$:
\[
G_\CC\ra \varphi_{\al*}(P_{|U_\al})\lra \varphi_{\al}(U_\al).\]
There is the commutative diagram:
\begin{equation}\label{eq:inducedbundle}
\begin{CD}
G_\CC@>>> G_\CC\\
@VVV @VVV \\
\varphi_{\be*}(P_{|U_\be})@>g_{\al\be*}>> \varphi_{\al*}(P_{|U_\al})\\
@VVV @VVV \\
\varphi_{\be}(U_\al\cap U_\al)@>g_{\al\be}>> \varphi_{\al}(U_\al\cap
U_\al).
\end{CD}
\end{equation}
Since the subgroup $\hat{\rm U}(n+1,1)$ acts transitively on
$S^{2n+1,1}$, we choose an element $h\in\hat{\rm U}(n+1,1)$ for
which $g=h\cdot g_{\al\be}\in {\rm PO}(2n+2,2)$ satisfies that
$gx=x$ for some point $x\in S^{2n+1,1}$. Then the differential map
$g_*:T_xS^{2n+1,1}\ra T_xS^{2n+1,1}$ satisfies that
$g_*\in G_\CC$.

Suppose that $H$ is a subgroup of $\mathrm {PO}(2n+2,2)$ containing
$\hat{\rm U}(n+1,1)$ which preserves the $G_\CC$-structure. As
above, note that $g\in H_x$. If $\tau:H_x\ra {\rm
Aut}(T_xS^{2n+1,1})$ is the tangential representation, then it
follows that
\[
\tau(H_x)\leq G_\CC\cong \RR^{2n}\rtimes (\mathrm{U}(n)\times
\RR^+).
\]
Since $\tau$ is injective for any connected compact subgroup of
$H_x$, this implies that a maximal compact subgroup $K'$ of $H_x$ is
isomorphic to ${\rm U}(n)$. Let $K$ be a maximal compact subgroup of
$H$ containing $K'$. By the Iwasawa-Levi decomposition,
$$K/K'\cong
H/H_x=S^{2n+1,1}=S^{1}\times S^{2n+1}/\ZZ_2,$$ $K$ must be
isomorphic
to ${\rm U}(n+1)\cdot \mathrm{U}(1)$.\\
On the other hand, $\displaystyle {\rm U}(n+1)\cdot {\rm U}(1)\leq
({\rm O}(2n+2)\cdot {\rm O}(2))$ is the maximal compact unitary
subgroup of $\hat{\mathrm{U}}(n+1,1)$. As $\hat{\rm U}(n+1,1)\leq
H$, we obtain that $\hat{\rm U}(n+1,1)=H$. In particular, $g=h\cdot
g_{\al\be}\in H_x\leq \hat{\rm U}(n+1,1)$. It follows that
$g_{\al\be}\in \hat{\rm U}(n+1,1)$.
Therefore the maximal collection of charts
$\{U_\al,\varphi_\al\}_{\al\in\Lambda}$ gives a uniformization with
respect to $(\hat {\rm U}(n+1,1),S^{2n+1,1})$.
\end{proof}

\begin{rremark}\label{rem:nilinfty}
 If $\hat\infty$ is the infinity point of $S^{2n+1,1}$ $($which
maps to the point at infinity $\{\infty\}$ of $S^{2n+1})$ $(${\rm
\cf} \eqref{inftypoint} of Section $\ref{sub:nocompactG})$, then it
is noted that the stabilizer $($up to conjugacy$)$ is
\[{\rm PO}(2n+2,2)_{\hat\infty}=
\RR^{2n+2}\rtimes ({\rm O}(2n+1,1)\times \RR^+).
\]
Note that the intersection $\hat{\rm U}(n+1,1)\cap{\rm
PO}(2n+2,2)_{\hat\infty}$ is $$\hat{\rm
U}(n+1,1)_{\hat\infty}=\mathcal N\rtimes ({\rm U}(n)\times \RR^+).$$
In fact, ${\rm O}(2n+1,1)$ contains the similarity subgroup
$\RR^{2n}\rtimes ({\rm O}(2n)\times\RR^+)$ so that
\[\mathcal N\rtimes {\rm U}(n)\subset
(\RR^{2n+2}\rtimes \RR^{2n})\rtimes {\rm O}(2n).\]

\end{rremark}
\vskip0.2cm

 The conformally
flat Lorentz geometry $({\rm PO}(2n+2,2),S^{2n+1,1})$ restricts a
subgeometry $(\hat{\rm U}(n+1,1),S^{2n+1,1})$. It is noted that the
full subgroup of ${\rm PO}(2n+2,2)$ preserving the $G$-structure on
$S^{2n+1,1}$ is $\hat{\rm U}(n+1,1)$.

\begin{ddefinition}\label{LCRLorentzG} The pair $(\hat{\rm
U}(n+1,1),S^{2n+1,1})$ is said to be \emph{conformally flat
Fefferman-Lorentz parabolic geometry}. A smooth $(2n+2)$-dimensional
manifold $M$ is a \emph{conformally flat Fefferman-Lorentz parabolic
manifold} if  $M$ is locally modelled on $\tilde S^{2n+1,1}$ with
local changes lying in ${\rm U}(n+1,1)^{\sim}$.\end{ddefinition}
 Here ${\rm U}(n+1,1)^{\sim}$ is a lift of
$\hat{\rm U}(n+1,1)$ to ${\rm PO}(2n+1,2)^{\sim}$ which has the
central group extension:

\begin{equation*}
1\ra \ZZ\ra {\rm U}(n+1,1)^{\sim}\stackrel{\hat Q_\ZZ}\lra \hat{\rm
U}(n+1,1)\ra 1.
\end{equation*}

We close this section by showing the following examples of
\emph{compact conformally flat Lorentz parabolic manifolds}.
\begin{itemize}
\item Lorentz  flat space forms which admit Lorentz parabolic structure but not
Fefferman-Lorentz parabolic structure.
\item Conformally flat Fefferman-Lorentz parabolic manifold which do not admit Fefferman-Lorentz
structure.
\item
Conformally flat Fefferman-Lorentz manifolds $S^1\times \mathcal
N^3/\Delta$ on which $S^1$ acts as lightlike isometries. (This is
shown in Section \ref{fl-sp} and Proposition \ref{Lorentznilflat}.)
\end{itemize}

\subsection{Examples of $4$-dimensional Lorentz flat parabolic
manifolds II}\label{examples2} Let $\mathcal N^3=\RR\times \CC$ be
the $3$-dimensional Heisenberg  group with group law:
\[(a,z)(b,w)=(a+b-{\rm Im}\,\bar zw,z+w).\]

Recall that Lorentz flat geometry $({\rm E}(2,1),\RR^3)$ where ${\rm
E}(2,1)=\RR^3\rtimes {\rm O}(2,1)$. Let ${\rm O}(2,1)_\infty$ be the
stabilizer at the point at infinity in $S^1=\partial \HH^2_\RR$. It
is isomorphic to ${\rm Sim}(\RR^1)= \RR\rtimes ({\rm O}(1)\times
\RR^+)$ as in \eqref{similaritygroup}. Put $z=x+{\mathrm i}t$,
$x,t\in \RR$. We define a continuous homomorphism:

\begin{equation}\label{flathomo}
\begin{split}
&\rho :\mathcal N^3\lra \RR^3\rtimes {\rm O}(2,1)_\infty,
\\
&\rho(\left(\begin{array}{c}
a\\
x\\
\end{array}\right))
= (\left(\begin{array}{c}
a\\
x\\
0\\
\end{array}\right),{\rm I}),\\
&\rho\left(\begin{array}{c}
{\mathrm i}t\\
\end{array}\right)=
(\left(\begin{array}{c}
-\mbox{\large$\frac{t^3}6$}\\
-\mbox{\large$\frac{t^2}2$}\\
t\\
\end{array}\right),
\left(\begin{array}{lcr}
1&t&-\mbox{\large$\frac{t^2}2$}\\
0&1&-t\\
0&0&1\\
\end{array}\right)).
\end{split}
\end{equation}
It is easy to see that $\rho$ is a simply transitive representation
of $\mathcal N^3$ onto the Lorentz flat space $\RR^3$. (Compare
\cite{kam5}.)

\begin{ppro}\label{eg1} There is a $4$-dimensional compact Lorentz flat space form $S^1\times
\mathcal N^3/\Delta$ which admits a Lorentz parabolic structure but
not admit Fefferman-Lorentz parabolic structure.\end{ppro}

\begin{proof}Taking $\RR$ as timelike parallel translations,
we extend the representation $\rho$ naturally to a simply transitive
$4$-dimensional representation:
\[\tilde\rho: \RR\times \mathcal N^3\lra
\RR\times \RR^3\rtimes {\rm O}(2,1)_\infty\subset  {\rm E}(3,1)\]

\begin{equation}\label{4dimflat}
\tilde \rho(\RR\times \mathcal N^3)=\RR^4.
\end{equation}Here note that ${\rm E}(3,1)
\RR^4\rtimes {\rm O}(3,1)\subset{\rm O}(4,2)_\infty$.
 If we choose a discrete uniform subgroup
$\Delta\subset \mathcal N^3$, then a compact aspherical manifold
$S^1\times \mathcal N^3/\Delta$ admits a (complete) flat Lorentz
structure such that
\[
S^1\times \mathcal N^3/\Delta\cong \RR^4/\tilde \rho(\ZZ\times
\Delta).\]We check that $S^1\times \mathcal N^3/\Delta$ cannot admit
a Fefferman-Lorentz parabolic structure. For this, if so, by
Proposition \ref{pro:flat-almost complex}, the group $\RR\times
\mathcal N^3$ is conjugate to a subgroup of ${\rm U}(2,1)$ up to an
element of ${\rm O}(4,2)$. Since $\RR\times \mathcal N^3$ is
nilpotent, it belongs to ${\rm U}(2,1)_\infty=S^1\cdot\mathcal
N^3\rtimes( {\rm U}(1)\times\RR^+)$ up to conjugate. This is
impossible because $S^1$ is lightlike.

\end{proof}

\subsection{Examples of conformally flat Fefferman-Lorentz parabolic
manifolds III}\label{examples1}
 We shall give compact conformally flat
Fefferman-Lorentz parabolic manifolds which are not equivalent to
the product of $S^1$ with spherical $CR$-manifold, \ie not a
Fefferman-Lorentz manifold. Consider the commutative diagram.
{\small
\begin{equation}\label{irrational-equiv}
\begin{CD}
\ZZ @. \ZZ@. \\
@VVV  @VVV@.\\
{\bf R}@>>>({\rm U}(n+1,1)^{\sim},\tilde S^{2n+1,1})
\searrow{}^{(\tilde P,\tilde P)}@.\\
@VVV  @V{(\hat Q_\ZZ,\hat Q)}VV  ({\rm PU}(n+1,1),S^{2n+1})\\
 S^1 @>>> (\hat{\rm
U}(n+1,1),S^{2n+1,1}) \nearrow{}_{(\hat P,P)}@.
\end{CD}
\end{equation}
}We start with a discrete subgroup $\Gamma\subset \hat{\rm
U}(n+1,1)$ such that $S^1\cap \Gamma=\ZZ_p$ for some integer $p$. If
we let $\pi={\hat Q_\ZZ}^{-1}(\Gamma)$, then there is the nontrivial
group extensions:
\begin{equation}\label{regulra-exact} {\small
\begin{CD}
1@>>> \frac 1p\ZZ@>>> \pi@>>> \tilde P(\Gamma)@>>> 1\\
 @.\cap @.\cap @.\cap \\
1@>>> {\bf R}@>>>{\rm U}(n+1,1)^{\sim}@>{\tilde P}>> {\rm PU}(n+1,1)@>>>1\\
\end{CD} }
\end{equation}
The group $\pi$ defines a cocycle $\displaystyle [f]\in H^2(\tilde
P(\Gamma),\frac 1p\ZZ)$. Suppose that $a$ is an irrational number.
Then
 $\displaystyle [a\cdot f]\in
H^2(\tilde P(\Gamma),{\bf R})$ which induces a group extension:
\[ 1\ra
\frac a{p}\ZZ \ra\pi(a)\lra \tilde P(\Gamma)\ra 1.\] Here $\pi(a)$
is viewed as the product $\displaystyle \frac a{p}\ZZ \times
P(\Gamma)$ with group law:
\[{\small
(\frac a{p}m, \al)(\frac a{p}\ell, \be)= (\frac a{p}(m+\ell)+ a\cdot
f(\al,\be),\al\be)}\ \ (\forall \al,\be\in \tilde P(\Gamma)).\]
(Refer to \cite{l-r2} and references therein for a construction of
group actions by group extensions.)

 As ${\bf R}$ is the center of
${\rm U}(n+1,1)^{\sim}$, it follows that
\begin{equation*}
(\frac a{p}m, \al)=(\frac a{p}m,1)(1, \al)\in {\bf R}\cdot{\rm
U}(n+1,1)^{\sim}={\rm U}(n+1,1)^{\sim}.
\end{equation*}This shows that
\begin{equation}\label{subgroup}
\pi(a)\subset {\rm U}(n+1,1)^{\sim}. \end{equation} As $\tilde
P(\Gamma)$ is discrete, so is $\pi(a)$ in ${\rm U}(n+1,1)^{\sim}$.
Let $L(\tilde P(\Gamma))$ be the limit set of $\tilde P(\Gamma)$ in
$S^{2n+1}$. Then it is known that $\tilde P(\Gamma)$ acts properly
discontinuously on the domain $\Omega=S^{2n+1}-L(\tilde P(\Gamma))$
(\cf \cite{ka-tsu},\cite{go}). If $\Omega\neq\emptyset$, then the
quotient $\Omega/\tilde P(\Gamma)$ is a spherical $CR$-orbifold.
Since $\displaystyle {\bf S}^1={\bf R}/\frac a{p}\ZZ$ is compact, it
is easy to see that $\pi(a)$ acts properly discontinuously on
$\tilde S^{2n+1,1}-{\tilde P}^{-1}(L(\tilde P(\Gamma)))$. Putting
\[
M(a)=\tilde S^{2n+1,1} -{\tilde P}^{-1}(L(\tilde P(\Gamma)))/\pi(a),
\]$M(a)$ is a smooth compact conformally flat Fefferman-Lorentz
 parabolic manifold which supports a fibration:
\[
{\bf S}^1 \ra M(a)\stackrel{\hat P}\lra \Omega/\tilde P(\Gamma).
\]
On the other hand, as $\displaystyle\hat Q_\ZZ(\pi(a))= \hat
Q_\ZZ(\frac a{p}\ZZ )\cdot \Gamma$ from \eqref{irrational-equiv},
the closure in $\hat{\rm U}(n+1,1)$ becomes
\begin{equation}\label{dense}
\overline {\hat Q_\ZZ(\pi(a))}=S^1\cdot \Gamma.
\end{equation}Whenever $a$ is irrational,
$M(a)$ cannot descend to a locally smooth orbifold modelled on
$(\hat{\rm U}(n+1,1),S^{2n+1,1})$. So $M(a)$ is not equivalent to
the product manifold. Hence we have

\begin{ppro}\label{pro:cfl-CR} Let $a$ be an irrational number.
There exists a compact $(2n+2)$-dimensional conformally flat
Fefferman-Lorentz parabolic manifold $M(a)$ which is a nontrivial
$S^1$-bundle over a spherical $CR$-manifold. Moreover, $M(a)$ is not
equivalent to the product manifold. \end{ppro}

For example, such $\pi(a)$ is obtained as follows. ${\rm PU}(n+1,1)$
has the subgroup ${\rm U}(n,1)={\rm P}({\rm U}(n,1)\times {\rm
U}(1))$ which acts transitively on
$S^{2n+1}-S^{2n-1}=V_{-1}^{2n,1}$. (See \eqref{cover-hyeprbolic}.)
Since the stabilizer at a point is isomorphic to ${\rm U}(n)$, there
exists a ${\rm U}(n,1)$-invariant Riemannian metric on
$V_{-1}^{2n,1}$. 
 If $\Gamma$ is a discrete cocompact subgroup of ${\rm U}(n,1)$,
 then $L(\Gamma)=L({\rm U}(n,1))=S^{2n-1}$. Chasing the diagram
\begin{equation}\label{discrete-diagram}
\begin{CD}
S^1@>>> \hat {\rm U}(n+1,1) @> \hat P>> {\rm PU}(n+1,1)\\
\bigcup @.        \bigcup @.   \bigcup @.\\
S^1   @>>> {\rm U}(n,1)\times {\rm U}(1)@>>> {\rm U}(n,1),
\end{CD}
\end{equation}
we start with $\Gamma \subset {\rm U}(n,1)\times \{1\}\subset {\rm
U}(n,1)\times {\rm U}(1)$. Then we get a Fefferman Lorentz manifolds
$\displaystyle M(a)={\bf R}\times V_{-1}^{2n+1}/\pi(a)$ where
$\tilde S^{2n+1,1}-\tilde S^{2n-1,1}={\bf R}\times V_{-1}^{2n+1}$.

Put $\pi=\pi(1) $ for $a=1$. In this case, $\hat Q_\ZZ(\pi)=\Gamma$.
The previous construction
 shows that
\begin{equation*}
\begin{split}
M(1)=\tilde S^{2n+1,1}-\tilde S^{2n-1,1}/\pi&=
S^{2n+1,1}-S^{2n-1,1}/\Gamma=
S^1\mathop{\times}_{\ZZ_2}^{}V_{-1}^{2n+1}/\Gamma.
\end{split}
\end{equation*}$M(1)$ is a conformally flat
Fefferman Lorentz manifold $M$. Varying $a$, we see that $M(a)$ is
nonequivalent with $M(1)$ as a Fefferman-Lorentz metric.
\\

\begin{rremark}\label{rem:compactpara}
We have also ${\rm O}(m+1)\times \RR^+$-structure on
$(m+2)$-manifolds as a parabolic structure. Similar to the proof of
Proposition \ref{pro:flat-almost complex}, we can show that
\begin{ppro}\label{pro:compact}
Let $M$ be a smooth $(m+2)$-manifold with an ${\rm O}(m+1)\times
\RR^+$-structure. If $M$ is conformally flat Lorentz such that ${\rm
O}(m+1)\leq {\rm O}(m+2)$ $(\emph{not maximal})$, then $M$ is
uniformized with respect to $(\mathrm{O}(2)\times \mathrm{O}(m+2),
S^1\times S^{m+1})$. In particular, if $M$ is compact, then $M$
covers $S^{m+1,1}$.
\end{ppro}
Let $M$ be a Riemannian manifold of dimension $m+1$. Then $S^1\times
M$ admits a natural Lorentz metric $\mathsf{g}$ for which $S^1$ acts
as timelike isometries. Even if $M$ is conformally flat, $S^1\times
M$ need not be a conformally flat Lorentz manifold. For example,
$S^1\times \HH^{m+1}_\RR$. However, $S^1\times \HH^{m+1}_\RR/\Gamma$
is covered by $({\rm P}({\rm O}(1,1)\times {\rm O}(m+1,1)),
\RR\times\HH^{m+1}_\RR)$ for which ${\rm P}({\rm O}(1,1)\times {\rm
O}(m+1,1))\leq {\rm PO}(m+2,2)$. So $S^1\times \HH^{m+1}_\RR/\Gamma$
is conformally flat Lorentz but $S^1$ $($or $\RR$ $)$ is not a group
of timelike isometries.
\end{rremark}

\section{Conformally flat Fefferman-Lorentz parabolic
geometry}\label{unitary}
 Recall that $\mathop{\Conf}_{\rm FLP}(M)$ is
the group of conformal transformations preserving the
Fefferman-Lorentz parabolic structure (\cf \eqref{fl-structure}). We
shall consider the representations of one-parameter subgroups
$H\leq{\Conf}_{\rm FLP}(M)$.
\subsection{One-parameter subgroups in $\hat{\rm
U}(n+1,1)$}\label{stab-basis}
 The following commutative diagrams are obtained.
{\footnotesize
\begin{equation}\label{groupgram}
\begin{CD}
\ZZ @=  \ZZ \\
@VVV       @VVV   \\
{\bf R} @>>> ({\rm U}(n+1,1)^{\sim},\tilde S^{2n+1,1})
@>(\tilde P,\tilde P)>>({\rm PU}(n+1,1),S^{2n+1})\\
@VVV  @V({\hat Q_\ZZ},\hat Q)VV     {||}\\
S^1 @>>>(\hat{\rm U}(n+1,1),S^{2n+1,1}) @>(\hat P,P)>> ({\rm PU}(n+1,1),S^{2n+1}).\\
\end{CD}
\end{equation}
} Here $\tilde S^{2n+1,1}=\RR\times S^{2n+1}$.

 {\footnotesize
\begin{equation}\label{groupgram1}
\begin{CD}
@. ({{\rm U}(n+1,1)^{\sim}}, \RR\times S^{2n+1})  \\
@. @V({Q_\ZZ},\tilde Q)VV  \searrow (\mbox{\tiny ${\hat Q_\ZZ},\hat
Q$})\quad \
\ \ \ \ \ \ \ \ \ \  \\
\ZZ_2 @>>> ({\rm U}(n+1,1),S^1\times S^{2n+1}) @>(Q_2,P_\RR)>> (\hat{\rm U}(n+1,1),S^{2n+1,1}).\\
\end{CD}
\end{equation}
}
By \eqref{good daigram}, there is the projection:
{\small
\begin{equation}\label{compositeP}
\begin{CD}
 P_{\CC}=P\circ P_{\RR}: V_0\, (\supset S^1\times S^{2n+1}) @>{P_\RR}>> S^{2n+1,1}@>P>> S^{2n+1}.
\end{CD}
\end{equation}
}As usual the following points $\{\infty,0\}$ are defined on the
conformal Riemannian sphere:
{\small
\begin{equation}\label{inf-north}\begin{split}
&\infty=P_\CC(f_1)=\left[\frac 1{\sqrt 2},0,\dots,0,\frac 1{\sqrt 2}\right]=
(0,\dots,0,1)\in S^{2n+1},\\
&0=P_\CC(f_{n+2})=\left[\frac 1{\sqrt 2},0,\dots,0,-\frac 1{\sqrt 2}\right]=
(0,\dots,0,-1)\in S^{2n+1}.
\end{split}\end{equation}
}
We put
\begin{equation*}\label{inftypoint}
\begin{split}
\hat\infty&=P_\RR(f_1)\in S^{2n+1,1},\ \, \ \hat 0=P_\RR(f_{n+2})\in
S^{2n+1,1}\\
\end{split}
\end{equation*}such that
\begin{equation*}
P(\hat\infty)=\infty, \ \ P(\hat 0)=0.
\end{equation*}

Suppose that $H$ is a one-parameter subgroup $\{\phi_t\}_{t\in\RR}$
of $\mathop{\Conf}_{\rm FLP}(M)$ and $\tilde
H=\{\tilde\phi_t\}_{t\in\RR}$ is its lift to ${\Conf}_{\rm
FLP}(\tilde M)$. 
 Let $\tilde\rho :\tilde
H\ra {\rm U}(n+1,1)^{\sim}$ be a homomorphism. For simplicity write
$\tilde \rho(\tilde \phi_t)=\tilde \rho(t)$ $(t\in \RR)$ and put
 \begin{equation}\label{holonomy}
\begin{split}
\rho&=Q_\ZZ\circ \tilde\rho \colon\tilde H\ra {\rm U}(n+1,1),\\
\hat\rho&=Q_2\circ \rho=\hat Q_\ZZ\circ \tilde\rho \colon\tilde H\ra
\hat{\rm U}(n+1,1).\\
\end{split}\end{equation}As
$P\colon {\rm U}(n+1,1)\to \mathrm {PU}(n+1,1)$ is the
projection, it follows from \eqref{groupgram}, \eqref{holonomy}
that $\displaystyle P\circ Q_\ZZ=\hat P\circ \hat Q_\ZZ=\tilde P$
for which
\begin{equation}\label{relation}
\begin{split}
\tilde P(\tilde \rho(t))&=P\rho(t).
\end{split}\end{equation}
Since $\tilde P\circ \tilde \rho(\tilde H)= P\circ\rho(\tilde H)$,
we put
\begin{equation}\label{cloG}
\mathsf {G}=\overline{P\circ\rho(\tilde H)} \leq{\rm PU}(n+1,1).
\end{equation}

We determine the connected closed subgroup
$\mathsf{G}$ by using the results of \cite{kam2}.
First recall that $\{e_1,\dots,e_{n+2}\}$ is the
standard complex basis of $\CC^{n+2}$ equipped with the Lorentz
Hermitian inner product $\langle\, ,\,\rangle$ (\cf \eqref{c-form}); $\langle
e_i,e_j\rangle=\delta_{ij}$ $(2\leq i,j\leq n+2)$, $\langle
e_{n+2},e_{n+2}\rangle=-1$. Setting $\displaystyle
f_1=e_1+e_{n+2}/\sqrt 2,\ \ f_{n+2}=e_1-e_{n+2}/\sqrt 2$
as before, the frame
$\displaystyle \{f_1,e_2,\dots,e_{n+1},f_{n+2}\}$ is the new basis such that
$$\langle f_1,f_1\rangle=\langle f_{n+2},f_{n+2}\rangle=0, \langle
f_1,f_{n+2}\rangle =\langle f_{n+2},f_{1}\rangle =1.$$

\subsubsection{Case I: $\mathsf{G}$ is
noncompact}\label{sub:nocompactG} It follows from \cite[\S 3]{kam2}
that $P(\rho(\tilde H))$ itself is closed. We may put
\begin{equation}\label{repreof G}
\mathsf{G}=P\circ \rho(\tilde H)=\{P\rho(t)\}_{t\in \RR}.
\end{equation}Moreover, $\mathsf{G}$ belongs to $\mathcal N\rtimes({\rm U}(n)\times \RR^+)=
{\mathrm PU}(n+1,1)_{\infty}$ up to conjugate. $($See Remark
\ref{rem:nilinfty}.$)$ Moreover, the explicit form of $\{P\rho(t)\}$
can be described with respect to the basis
$\{f_1,e_2,\cdots,e_{n+1},f_{n+2}\}$. (Compare \cite{kam3}.) It has
the following form
\begin{equation}\label{nil} P\rho(t)=\left[\begin{array}{ccc}
 1 &  \mbox{\Large $0$}   & t\mbox{\boldmath$i$}\\
\mbox{\Large $0$} & A_t   & \mbox{\Large $0$} \\
0 & \mbox{\Large $0$} & 1
\end{array}\right]
\end{equation}where $\displaystyle A_t=(e^{\mathrm{i}ta_1},\dots,
e^{\mathrm{i}ta_n})\in T^n\leq{\rm U}(n)$.
\begin{equation}\label{c-nil}
P\rho(t)=\left[\begin{array}{cccc}
1 &                t &\mbox{\Large $0$} & t^2/2+1\\
0 &                1 & \mbox{\Large $0$}& t \\
\mbox{\Large $0$}&\mbox{\Large $0$} & B_t              & \mbox{\Large $0$} \\
0 & 0            & \mbox{\Large $0$} & 1
\end{array}\right]
\end{equation}where $\displaystyle B_t=(e^{\mathrm{i}tb_1},\dots,
e^{\mathrm{i}tb_{n-1}})\in T^{n-1}\leq{\rm U}(n-1)$.

\begin{equation}\label{com}
P\rho(t)=\left[\begin{array}{ccc}
 e^t &  \mbox{\Large $0$}   & 0\\
 \mbox{\Large $0$}& A_t   & \mbox{\Large $0$} \\
0 &\mbox{\Large $0$} & e^{-t}
\end{array}\right]
\end{equation}where $\displaystyle A_t=(e^{\mathrm{i}ta_1},\dots,
e^{\mathrm{i}ta_n})\in T^n\leq{\rm U}(n)$.
\\

Let $C_t$ be the matrix accordingly as whether  $[C_t]$ is
\eqref{nil}, \eqref{c-nil} or \eqref{com}. Noting that the center of
${\rm U}(n+1,1)$ is $S^1=\{e^{\mathrm{i}t}\}$, the holonomy map
$\rho:\tilde H\ra {\rm U}(n+1,1)$ has the following form:

\begin{equation}\label{nil-lift}
\rho(t)=\left\{\begin{array}{lr}
C_t & \ \ \ {\rm (i)},\\
 e^{\mathrm{i}t}\cdot C_t &{\rm (ii)} .
\end{array}\right.
\end{equation}
For \eqref{nil}, \eqref{c-nil}, $\mathsf{G}$ has the unique fixed
point $\{\infty\}$ in $S^{2n+1}$.
As $P_\CC(f_1)=\infty$ (\cf
\eqref{inf-north}), $\rho(t)f_1=\lambda\cdot f_1$ for some
$\lambda\in \CC^*$. If $\rho(t)=C_t$ for \eqref{nil-lift}, then
$\rho(t)f_1=f_1$ so that $\hat\rho(t)\hat\infty=\hat\infty$ by
\eqref{groupgram1}. Hence
\begin{equation}\label{CaseA}
\hat\rho(\tilde H)\ \mbox{has the fixed point set}\ \{S^1\cdot
\hat\infty\}\ \mbox{in}\ S^{2n+1,1}.
\end{equation}
 For \eqref{com}, $\mathsf{G}$ has two fixed
points $\{0, \infty\}$ in $S^{2n+1}$. If $\rho(t)=C_t$, then
$\rho(t)f_1=e^t\cdot f_1$, $\rho(t)f_{n+2}=e^{-t}\cdot f_{n+2}$.
Since $P_\RR(s\cdot v)=P_\RR(v)$ for $\forall s\in
\RR^*,  v\in V_0$, it follows that
$\hat\rho(t)\hat\infty=\hat\infty$ and $\hat\rho(t)\hat 0=\hat 0$ in
$S^{2n+1,1}$. Similarly as above,
\begin{equation}\label{CaseB}
\hat\rho(\tilde H)\ \mbox{has the fixed point set}\ \{S^1\cdot\hat
0, S^1\cdot\hat\infty\}\ \mbox{in}\ S^{2n+1,1}.
\end{equation}

\subsubsection{Case II: $\mathsf{G}$ is compact}\label{sub:compactG}
 Using \eqref{relation},
\[
P(\rho(t))=
(e^{\mbox{\boldmath$i$}ta_1},\dots,e^{\mbox{\boldmath$i$}ta_k},1,\dots,1)\in
T^{n+1}\]for some nonzero numbers $a_1,\dots, a_k$. Put
$E_t=(e^{\mbox{\boldmath$i$}ta_1},\dots,e^{\mbox{\boldmath$i$}ta_k},1,\dots,1)$.
We may assume that the g.c.m of all $a_i$ is $1$ (up to scale factor
of parameter $t$). Then $\rho(t)$ has one of the following forms:
\begin{equation}\label{compact-lift1}
\rho(t)=\left\{\begin{array}{lr} E_t \in
T^{n+1}\cdot S^1 \ \ \ \ \ \ \ \ \ \ {\rm (i)},\\
E_t\cdot e^{\mathrm{i}t}\in T^{n+1}\cdot S^1 \ \ \ \ {\rm (ii)}.
\end{array}\right.
\end{equation}
\begin{ppro}\label{G-noncompact}\
Let $(M,g)$ be a conformally flat Fefferman-Lorentz parabolic
 manifold which admits a one-parameter subgroup
$H\leq\mathop{\Conf}_{\rm FLP}(M)$ acting without fixed points on
$M$. Suppose that
\begin{equation}\label{devpairH} (\tilde\rho,\dev):(\tilde H,\tilde
M)\ra ({\rm U}(n+1,1)^{\sim},\tilde S^{2n+1,1})
\end{equation}is the developing pair.  If $\mathsf{G}\neq \{1\}$, then either one
of {\bf { Case A},
{Case B}, { Case C}} or {\bf { Case D}} holds:
\par\ \par

\begin{itemize}
\item[{\bf {Case A.}}] The action of $\rho(\tilde H)$ is
 of type ${\rm (i)}$ of \eqref{nil-lift}.
\begin{itemize}
\item[(1)]\ When $C_t$ has the form of either \eqref{nil} or \eqref{c-nil},

\begin{equation}\label{fixout}
\hat Q(\mathop{\dev}(\tilde M)) \subset
S^{2n+1,1}-S^1\cdot{\hat\infty}=S^1\times \mathcal N,
\end{equation}where
the centralizer $\mathcal C(\hat\rho(\tilde H))$ of
$\hat\rho(\tilde H)$ in $\hat {\rm U}(n+1,1)$ is contained in
$S^1\times (\mathcal N\rtimes {\rm U}(n))$.
\item[(2)]\ When $C_t$ has the form \eqref{com},
\begin{equation}\label{fixout1}
 \hat Q(\mathop{\dev}(\tilde M)) \subset
S^{2n+1,1}-S^1\cdot\{\hat 0,\hat\infty\}=(S^{2n}\times \RR^+)\times
S^1,\end{equation} where $\mathcal C(\hat\rho(\tilde H))$ is
contained in $S^1\times ({\rm U}(n)\times \RR^+)$.
\end{itemize}

\item[{\bf {Case B.}}] The action of $\rho(\tilde H)$ is of type
 ${\rm (ii)}$ of \eqref{nil-lift}. Then $\hat\rho(\tilde H)$ has no fixed point set
 on $S^{2n+1,1}$.
In this case, $\mathcal C(\hat\rho(\tilde H))$ is either contained
in $S^1\times (\mathcal N\rtimes {\rm U}(n))$ or $S^1\times ({\rm U}(n)\times \RR^+)$.

\item[{\bf {Case C.}}] The action of
$\hat\rho(\tilde H)$ is of type ${\rm (i)}$ of \eqref{compact-lift1}.
\begin{equation}\label{compactimage1}
 \hat Q(\mathop{\dev}(\tilde M)) \subset
S^{2n+1,1}-S^1\cdot S^{2(n-k)+1}
\end{equation}on which the subgroup
$\displaystyle S^1\times
({\rm U}(n-k+1,1)\widehat{\times} {\rm U}(k))$ acts transitively with compact stabilizer.
The centralizer of
$\hat\rho(\tilde H)$ in $\hat{\rm U}(n+1,1)$ is $\displaystyle
S^1\times({\rm U}(n-k+1,1)\widehat{\times} T^k)$ where
$\hat\rho(\tilde H)\subset T^k$.

\item[{\bf {Case D.}}]
The action of $\hat\rho(\tilde H)$ is of type ${\rm (ii)}$ of
\eqref{compact-lift1}. Then $\hat\rho(\tilde H)$
has no fixed point set  on $S^{2n+1,1}$ and
its centralizer in  $\hat{\rm U}(n+1,1)$ is
$\displaystyle S^1\times({\rm U}(n-k+1,1)\widehat{\times} T^k)$
 where $S^1\times\hat\rho(\tilde H)\subset T^k$.
\end{itemize}

\end{ppro}

\begin{proof}
By the hypothesis, $H$ has no fixed point so does $\tilde H$ on
$\tilde M$. Since $\mathop{\dev}$ is an immersion, the image
$\mathop{\dev}(\tilde M)$ misses the fixed point set of
$\tilde\rho(\tilde H)$ in $\tilde S^{2n+1,1}$. Recall that there is
the covering space from \eqref{groupgram}:
\[
\begin{CD}
 \ZZ\ra ({\rm  U}(n+1,1)^{\sim}, \tilde
S^{2n+1,1})@>(\hat Q_\ZZ,\hat Q)>> (\hat{\rm U}(n+1,1),S^{2n+1,1}).
\end{CD}\]
Noting that $\tilde\rho(\tilde H)$ is connected, the image $\hat
Q(\mathop{\dev}(\tilde M))$ also misses the fixed point set of $\hat
Q_\ZZ(\tilde\rho(\tilde H))=\hat\rho(\tilde H)$ of $\hat {\rm
U}(n+1,1)$ (\cf \eqref{holonomy}). Then $(1)$, $(2)$ of {{\bf Case
A}} follow from \eqref{CaseA} and \eqref{CaseB} respectively and
{\bf {Case B}} follows easily because the center $S^1$ of $\hat{\rm
U}(n+1,1)$ acts freely on $S^{2n+1,1}$.

For the case (i) of \eqref{compact-lift1}, the fixed point set of
$P(\rho(\tilde H))$ is
$$S^{2(n-k)+1}=\{(0,\dots,0,z_{k+1},\dots,z_{n+1})\in S^{2n+1}\}$$ in which
 the subgroup of ${\rm PU}(n+1,1)$ preserving
$S^{2n+1}-S^{2(n-k)+1}$ is ${\rm P}({\rm U}(n-k+1,1)\times {\rm
U}(k))$. Since $\hat Q(\mathop{\dev}(\tilde M))$ misses the fixed
point set of $\hat \rho(\tilde H)$,  {\bf Case C} follows that
\begin{equation*}
\begin{split}
&\hat\rho(t)= (e^{\mathrm{i}ta_1},\dots,e^{\mathrm{i}ta_k}, 1,
\dots,1)\in T^k\leq S^1\cdot T^{n+1},\\
&\hat Q(\mathop{\dev}(\tilde M)) \subset S^{2n+1,1}-S^1\cdot
S^{2(n-k)+1}. \end{split}\end{equation*}
 For the case (ii) of \eqref{compact-lift1}, it follows that
\[\hat\rho(t)
=(e^{\mathrm{i}t(a_1+1)},\dots,e^{\mathrm{i}t(a_k+1)},
e^{\mathrm{i}t}, \dots,e^{\mathrm{i}t})\in S^1\cdot T^{n+1}.\]
Similarly as {\bf Case B}, $$\mathcal C(\hat \rho(\tilde H))\leq
S^1\cdot {\rm U}(n-k+1,1)\widehat{\times} T^k)$$ where
$\hat\rho(\tilde H)\subset S^1\cdot T^k$ which shows {\bf Case D}.

\end{proof}

\noindent Denote by $\mathrm{Isom}_{\rm FLP}(M)$ the group of
isometries preserving the Fefferman-Lorentz parabolic structure such
that $\mathrm{Isom}_{\rm FLP}(M)\leq \mathop{\Conf}_{\rm FLP}(M)$.

\begin{ppro}\label{1G-noncompact}\
Let $(M,g)$ be a conformally flat Fefferman-Lorentz parabolic
 manifold admitting a $1$-parameter subgroup
$H\leq\mathrm{Isom}_{\rm FLP}(M)$. If $\mathsf{G}=\{1\}$, then the
following hold.
\begin{itemize}
\item[(i)] The lift $\tilde H$
acts properly and freely on $\tilde M$
as lightlike isometries.

\item[(ii)] $\tilde M/\tilde H$ is a simply connected
spherical $CR$-manifold on which the quotient group $\mathcal
C(\tilde H)/\tilde H$ acts as
$CR$-transformations. 
Here $\mathcal C(\tilde H)$ is the centralizer
of $\tilde H$ in $\mathop{\Conf}_{\rm FLP}(\tilde M)$.
\item[(iii)] The conformal developing pair
$(\tilde \rho,\mathop{\dev})$
for $M$ induces a $CR$-developing pair:
\begin{equation}\label{CRdev}
(\hat \rho,\hat{\dev})\colon ({\mathcal C}(\tilde H)/\tilde H,\tilde
M/\tilde H)\to ({\mathrm {PU}}(n+1,1),S^{2n+1}).
\end{equation}
\end{itemize}

\end{ppro}

\begin{proof}
Suppose that $\mathsf{G}=
\{1\}$.
Since $\mathsf{G}=
\overline{\tilde P\circ \tilde \rho(\tilde H)}$ from
\eqref{relation}, \eqref{cloG}, it follows that
$\tilde \rho(\tilde H)=\RR\leq {\mathrm{U}}(n+1,1)^{\sim}$
which is lightlike with respect to $g^0$ where
\begin{equation}\label{form}
g^0=\sigma_0\odot \tilde P^*\om_0+ d\om_0(J_0\tilde P_*-,\tilde
P_*-)
\end{equation}
is the standard Lorentz metric on
$\tilde S^{2n+1,1}=\RR\times S^{2n+1}$ induced from \eqref{f-lorentz}.
In particular, $\tilde \rho\colon\tilde H\to\RR$ is an
isomorphism. As $\RR$ acts properly on $\tilde
S^{2n+1,1}=\RR\times S^{2n+1}$, $\tilde H$ acts properly on $\tilde
M$.
On the other hand, there
exists a function $u$ such that
\begin{equation}\label{conf-u}
{\dev}^*g^0=u\cdot g.
\end{equation}
Let $\mathcal H$ be the vector field induced by $\tilde H$ on $\tilde M$.
As $\mathcal S$ is the vector field induced by
$\RR$, we have that $\mathop{\dev}_*(\mathcal H )=\mathcal
S$. Since $u\cdot
g(\mathcal H,\mathcal H)= g^0(\mathcal S,\mathcal S)=0$,
noting the hypothesis that
$\tilde H\leq\mathrm{Iso}_{\rm FLP}(\tilde M)$,
$\tilde H$ acts as lightlike isometries.
This shows (i).

There is the commutative diagram:
\begin{equation}\label{downCR}\begin{CD}
\tilde M@>\dev>> \tilde S^{2n+1,1}\\
@V{P}VV  @V{\tilde P}VV \\
\tilde M/\tilde H @>\hat \dev >> S^{2n+1}.
\end{CD}
\end{equation}

We put
\begin{equation}\label{contact}
P^*\om (X)=g(\mathcal H, X)\ \, (\forall\, X\in T\tilde M).
\end{equation}Then $\om$ is a well-defined $1$-form on $\tilde M/\tilde H$
by the fact that $g(\mathcal H,\mathcal H)=0$.

Note from
\eqref{sigmaform}, \eqref{leescomp} that
\begin{equation}\label{stanform}
\sigma_0=\frac 1{n+2}\left(dt+\mathrm{i}\tilde
P^*\pi^*\varphi_\al^{\al}-\frac {nc}4\cdot \tilde P^*\om_0\right).
\end{equation}For $\tilde \rho(\tilde H)=\RR$, it follows that
$\tilde \rho(h)^*\sigma_0=\sigma_0$ $(\forall
h\in \tilde H)$.
This implies that $\tilde \rho(h)^*g^0=g^0$. Applying $\mathop{\dev}^*$ to this,
$h^*\mathop{\dev}^*g^0=\mathop{\dev}^*g^0
=u\cdot g$ by \eqref{conf-u}.
As $h^*\mathop{\dev}^*g^0=h^*(u\cdot g)=
h^*u\cdot g$,
it follows
$h^*u=u$ and so $u$ factors through a map $\hat u\colon \tilde
M/\tilde H\to \RR^+$ such that
\begin{equation}\label{factors}
P^*\hat u=u.
\end{equation}

Moreover $\displaystyle
\sigma_0(\mathcal S)=\frac1{n+2}$ from
\eqref{dtvalue} and \eqref{stanform}, so \eqref{form} implies that
$$g^0(\mathcal S,-)=\frac 1{n+2}\tilde P^*\om_0(-).$$
Using \eqref{factors}, the equation
$\mathop{\dev}^*g^0=u\cdot g$ yields that
\begin{equation*}
\begin{split}
P^*\hat u\cdot P^*\om(X)&=u\cdot g(\mathcal H,X)
        =g^0(\mathcal S, \dev_*X)\\
   &=\frac 1{n+2}\tilde P^*\om_0(\dev_*X)=\frac 1{n+2} \dev^*\tilde P^*\om_0(X)\\
      &=\frac 1{n+2}P^*\hat{\dev^*}\om_0(X),
\end{split}\end{equation*}hence
\begin{equation}\label{crequiv}
(n+2)\hat u\cdot \om=\hat{\dev^*}\om_0.
\end{equation}
As $\hat\dev\colon\tilde M/\tilde H\to S^{2n+1}$ is an immersion,
$\hat\dev_*\colon \mathrm{ker}\, \om\to \mathrm{ker}\, \om_0$ is an
isomorphism. Define $\hat J$ on $\mathrm{ker}\, \om$ to be
\begin{equation}\label{complexCR}
\hat\dev_*(\hat JX)=J_0\hat\dev_*(X). \end{equation} If we note that
$J_0$ is a complex structure on $\mathrm{ker}\, \om_0$, $\hat J$
turns out to be a complex structure on $\mathrm{ker}\, \om$. Hence
$(\mathrm{ker}\, \om,\hat J)$ gives a $CR$-structure on $\tilde
M/\tilde H$ for which $\hat{\mathop{\dev}}$ is a $CR$-immersion.

Let $\mathcal C(\tilde H)$ be the centralizer
of $\tilde H$ in $\mathop{\Conf}_{\rm FLP}(\tilde M)$.
For $s\in {\mathcal C}(\tilde H)$ with
$\hat s\in {\mathcal C}(\tilde H)/\tilde H$,
there is a positive function $v$ on $\tilde M$
such that $s^*g=v\cdot g$. Noting that
$\tilde H\leq\mathrm{Iso}_{\rm FLP}(\tilde M)$,
we can check that $\tilde h^*v=v$
$(\forall h\in \tilde H)$, \ie there exists a function $\hat v$ on
$\tilde M/\tilde H$ such that $P^*\hat v=v$.
Then it is easy to see that
$\hat s^*\om=\hat v\cdot \om$ on $\tilde M/\tilde H$.
Using \eqref{complexCR},
it follows that $\hat s_*\circ\hat J=\hat J\circ \hat s_*$
on $\mathrm{ker}\, \om$.
Hence the group ${\mathcal C}(\tilde H)/\tilde H$
 preserves the $CR$-structure $(\mathrm{ker}\, \om,\hat J)$ on
$\tilde M/\tilde H$. This shows (ii).

As $\tilde \rho(\tilde H)=\RR$ is the center
of ${\mathrm{U}}(n+1,1)^{\sim}$, $\rho \colon
\mathrm{Conf}_{\rm FLP}(\tilde M)\to \mathrm{U}(2n+1,1)^{\sim}$
induces a homomorphism $\hat\rho\colon {\mathcal C}(\tilde
H)/\tilde H\to{\mathrm {PU}}(n+1,1)$.
Using the above commutative diagram, it follows that
\[
\hat{\dev}(\hat \al\hat x)=\hat\rho(\hat\al)
\hat{\dev}(\hat x)\ \, (\forall\, \hat\al\in {\mathcal C}(\tilde H)/\tilde H,
\hat x\in\tilde M/\tilde H).
\]
Hence (iii) is proved.

\end{proof}

\begin{rremark}\label{g-invariant}
A Fefferman-Lorentz manifold $M=S^1\times N$ is obviously an example
satisfying the hypothesis of Proposition $\ref{1G-noncompact}$.
\end{rremark}

\section{Coincidence of curvature flatness}\label{Fefferman-Lorentz-CR manifolds}
In this section we shall prove the equivalence between conformally
flatness of Fefferman-Lorentz manifolds and (spherical) flatness of
underlying $CR$-manifolds. (Compare Theorem
\ref{equivalence-flatness}.)

\subsection{Causal vector fields on $S^{2n+1,1}$}\label{conf-metric}
Let $H\leq \mathrm{Conf}_{FLP}(M)$ be as before and $\rho:\tilde H\to
\mathrm {U}(n+1,1)$ the representation. Denote by $\xi$ the vector
field on $S^1\times S^{2n+1}$ induced by the orbit $\rho(\tilde
H)\cdot z$ for some $z\in S^1\times S^{2n+1}$. Recall from
\eqref{groupgram} that
\begin{equation}\label{quotientvector}
P_{\RR}(\rho(t)(\lambda\cdot z))=Q_{2}(\rho(t))\cdot P_\RR
(\lambda\cdot z) =\hat \rho(t)\cdot P_\RR (z)
\end{equation}so we put the vector field $\hat \xi$ on
$S^{2n+1,1}$ by
\begin{equation}\label{eq:causalfield}
(P_{\RR})_*(\xi_{z})=\hat\xi_{P_{\RR}(z)}.
\end{equation} Note that $P_*(\hat \xi)=P_{\CC*}(\xi_{z})$
is a vector field on $S^{2n+1}$ by using
$P_{\CC}=P\circ P_{\RR}$ (\cf \eqref{compositeP}).
From \eqref{f-lorentz}, let
\[g^0(X,Y)=(\sigma_0\odot
{P_\CC}^*\omega_0)(X,Y)+d\om_0(J_0{P_\CC}_*X,{P_\CC}_*Y)\] be the
standard Lorentz metric on $S^1{\times}S^{2n+1}$ $(X,Y\in
T(S^1\times S^{2n+1})$).

Using the classification of one-parameter subgroups $\rho(\tilde
H)=\{\rho(t)\}\leq{\rm U}(n+1,1)$ of Section \ref{stab-basis},
 we examine the causality of the
vector field $\xi$ induced by $H$.\\

\subsection{Nilpotent group case}\label{nipotentcase}
(Compare Section \ref{sub:nocompactG},\eqref{nil-lift}.)

As $g_\mathcal N=u\cdot g^0$ on $S^1\times \mathcal N$, using
$g_\mathcal N$ instead of $g^0$, it suffices to check the causality
of $\xi$. Note that $\displaystyle\sigma_\mathcal N=\frac 1{n+2}dt$.

For the vector field $\xi$ restricted to $S^1\times \mathcal
N\subset S^1\times S^{2n+1}$, let $P_*(\xi)$ be the vector field on
$\mathcal N$ which is induced by the one-parameter subgroup
$P(H)\leq \mathrm{PU}(n+1,1)$ as above. It is not necessarily a
characteristic vector field except for the case $(i)$ of
\eqref{nil-lift}, but note that $\om_\mathcal N(P_{*}(\xi))\neq 0$
on $\mathcal N$. (Compare \cite{kam2}.) In fact, for the cases
\eqref{nil}, \eqref{c-nil}, \eqref{com} respectively,

\begin{equation*}\begin{split}
(1)\ &\ P_{*}(\xi)=\frac d{dt}
+\mathop{{\sum}}_{j=1}^{k}a_j(x_j\frac d{dy_j}-y_j\frac {d}{dx_j}),\\
&\omega_\mathcal N(P_{*}(\xi))=1+(a_1|z_1|^2+\cdots+a_k|z_k|^2).\\
(2)\ &\ P_{*}(\xi)=-y_1\frac d{dt}+\frac d{dx_1}+
\mathop{{\sum}}_{j=2}^{n}b_j(x_j\frac d{dy_j}-y_j\frac
{d}{dx_j}),\\
&\omega_\mathcal N(P_{*}(\xi))
=-2y_1+(b_2|z_2|^2+\cdots+b_{n}|z_n|^2).\\
(3)\ &\ P_{*}(\xi)=2t\frac d{dt}
+\mathop{{\sum}}_{j=1}^{n}((x_j-a_jy_j)\frac d{dx_j}+
(y_j+a_jx_j)\frac {d}{dy_j}),\\
&\omega_\mathcal
N(P_{*}(\xi))=2t+(a_1|z_1|^2+\cdots+a_n|z_n|^2).\\
\end{split}\end{equation*}

As above, $\displaystyle\sigma_\mathcal N(\xi)=\frac
1{n+2}dt(\xi)=\frac 1{n+2}\delta$ where $\delta=0,\, 1$ according to
the case (i) or (ii) of  \eqref{nil-lift} respectively. We obtain
that
\[
g_\mathcal N(\xi,\xi)=\frac {2\delta}{n+2}\om_\mathcal
N(P_{*}(\xi))+d\om_\mathcal N(JP_{*}(\xi),P_{*}(\xi)).\] Moreover,
it follows from \eqref{nil-con} that
\begin{equation}\label{metric}
d\omega_{\mathcal N}(JP_{*}X,P_{*}Y)=
\mathop{{\sum}}_{j=1}^{n}(dx_j^2+dy_j^2)(P_{*}X,P_{*}Y).
\end{equation} Calculating $d\om_\mathcal N(JP_{*}(\xi),P_{*}(\xi))$
for the above $P_*(\xi)$ respectively, we see that $g_\mathcal
N(\xi,\xi)=0$ if and only if $\delta=0$ and $P_{*}(\xi)$ is the
characteristic vector field
for $\om_\mathcal N$ in (1). As a consequence, we obtain that\\

\par\noindent {\bf Causality (1).}\ \, Suppose that
$\mathsf{G}$ is noncompact. Let $\xi$ be a lightlike conformal
vector field on $S^{2n+1,1}$ induced by $\mathsf{G}$. Then $\xi$ is
a nonzero lightlike vector field of $g_\mathcal N$ on $S^1\times
\mathcal N$ if and only if $\tilde H=\RR$, $\displaystyle
P_{*}(\xi)=\frac d{dt}$ is the
characteristic vector field for $\om_\mathcal N$.\\

\subsection{Compact torus case}\label{compactcase}
(\cf Section \ref{sub:compactG}, \eqref{compact-lift1}).
 It is possible to calculate $g^0$ by
making use of $\sigma_0$, however it is difficult to see
$\om_\al^\al$. So we consider a different approach. Let $\RR^*\ra
V_0\stackrel{P_{\RR}}\lra S^{2n+1,1}$ be the projection for which
$\hat g^0$ is the standard Lorentz metric on $S^{2n+1,1}$ with
$g^0=P^*_\RR \hat g^0$.  If we choose $\displaystyle c=\frac 2n$ in
\eqref{sigmasimple}, then
\begin{equation}\label{xixix}
\hat g^0(\xi,\xi)=-\frac 1{n+2}\ \, \mbox{on}\ \, S^{2n+1,1}.
\end{equation}

Recall from \eqref{c-form} that
\[V \sb 0=\{z=(z_1,\dots,z_{n+2})\in
\CC^{n+2}-\{0\}\ |\ \langle z,z\rangle=0\}\] in which
$\displaystyle\langle z,w\rangle =\bar z_1w_1+\cdots+\bar
z_{n+1}w_{n+1}-\bar z_{n+2} w_{n+2}\ \ \mbox{on}\ \CC^{n+2}$.
For an
arbitrary point $P_\RR(\bar p)=[\bar p]\in S^{2n+1,1}$ (\cf
\eqref{good daigram}), choose a point $\bar q\in V_0$ such that
\begin{equation}\label{productpq}
\begin{split}
\langle \bar p,\bar q\rangle= r\ \,(\exists\, r\in \RR-\{0\}).
\end{split}\end{equation}
Note that $\langle \bar p,\bar p\rangle=\langle \bar q,\bar
q\rangle=0$. From \eqref{good daigram}, we have the decomposition:
\begin{equation}\label{light}
\begin{split}
T_{\bar p}\RR^*=\RR\bar p &\ra T_{\bar p}V_0=\CC\bar
p+\mbox{\boldmath$i$}\RR\bar q+W_0\\
& \stackrel{P_{\RR*}}\lra T_{[\bar p]}S^{2n+1,1}=
\CC\bar p+\mbox{\boldmath$i$}\RR\bar q+W_0/\RR\bar p\\
&\quad \ \ \ \ \ \ \ \ \ \ \ \ \ \ \ \ \  \approx
\mbox{\boldmath$i$}\RR\bar p+\mbox{\boldmath$i$}\RR\bar q+W_0
\end{split}
\end{equation}where $T_{\bar p}V_0=\{v\in \CC^{n+2}\ |\
{\rm Re}\langle \bar p,v\rangle=0\}$ and $W_0={\langle \bar p,\bar
q\rangle}^{\perp}$. The decomposition is independent of the choice
of $\bar q$ with respect to $\langle \bar p,\bar q \rangle\neq 0$.
Note that ${\rm Re}\langle\ ,\ \rangle$ is the metric on $V_0$ of
dimension $2n+3$.

\begin{ppro}\label{anothermetric}
When $\displaystyle c=\frac 2n$,
\[ \hat g^0(P_{\RR*}X,P_{\RR*}Y)=
{\rm Re}\langle X,Y\rangle\ \ (X,Y\in T_{\bar p}V_0).\]
\end{ppro}

\begin{proof}
For an arbitrary point $x\in S^{2n+1,1}$, we can choose
\[ \bar p=(a_1,\dots,a_{n+1},\frac 1{\sqrt{n+2}}z)\ \, (\exists
\,z\in S^1)\] such that $P_\RR(\bar p)=[\bar p]=x$. Then $\xi_x$ is
induced by the $S^1$-orbit at $\bar p$:
 \[c(\theta)=(a_1,\dots,a_{n+1},\frac 1{\sqrt{n+2}}z\cdot e^{-\mathrm{i}\theta})\in V_0.\]
Since $\displaystyle \dot{c}(0)=(0,\cdots,0,-\frac
1{\sqrt{n+2}}z\mathrm{i})$, it follows that
$P_{\RR*}(\dot{c}(0))=\xi_x$. Similarly, $\mathcal S_x$ is induced
by the $S^1$-orbit at $\bar p$:
 \[s(\theta)=(a_1\cdot e^{\mathrm{i}\theta},\dots,a_{n+1}\cdot e^{\mathrm{i}\theta},
 \frac 1{\sqrt{n+2}}z\cdot e^{\mathrm{i}\theta})\in V_0.\]
 It follows that $\displaystyle \dot{s}(0)=(a_1\mathrm{i},\dots,a_{n+1}\mathrm{i},
 \frac 1{\sqrt{n+2}}z\mathrm{i})$
 for which $P_{\RR*}(\dot{s}(0))=\mathcal S_x$.
 Then we check that
 \begin{equation*}
 \begin{split}
{\rm Re}\langle \dot{c}(0),\dot{c}(0)\rangle&=-\frac 1{n+2},\\
{\rm Re}\langle \dot{s}(0),\dot{c}(0)\rangle&=\frac 1{n+2},\\
{\rm Re}\langle \dot{s}(0),\dot{s}(0)\rangle&=\langle \bar p,\bar p\rangle=0.\\
\end{split}
\end{equation*}
Thus \eqref{xixix} or \eqref{pairofLo} respectively shows that
\begin{equation*}\begin{split}
\hat g^0(\xi,\xi)&=\hat
g^0(P_{\RR*}(\dot{c}(0)),P_{\RR*}(\dot{c}(0)))= {\rm Re}\langle
\dot{c}(0),\dot{c}(0)\rangle.\\
\hat g^0(\mathcal S,\xi)&=\hat g^0(P_{\RR*}(\dot{s}(0)),
P_{\RR*}(\dot{c}(0)))={\rm Re}\langle
\dot{s}(0),\dot{c}(0)\rangle.\\
\end{split}\end{equation*}

Since it is easy to see that $\mathop{\ker}\,\om_0=P_{\RR*}(W_0)$
from \eqref{sphericalform}, we have that
\begin{equation*}
\begin{split}
d\om_0(JP_{\RR*}X,P_{\RR*}Y)&= |dz_1|^2+\cdots+|dz_{n+1}|^2 (X,Y)\\
&= {\rm Re}\langle X,Y\rangle\ \, (\forall X,Y\in W_0).
\end{split}\end{equation*}
Hence we obtain that
\[
\hat g^0(P_{\RR*}X,P_{\RR*}Y)={\rm Re}\langle X,Y\rangle \ \,
(X,Y\in T(V_0)).\]

\end{proof}

\begin{llemma}\label{lem:C-com}
If the one-parameter group has the form $(${\rm \cf}\ {\bf Case
C}$)$
$$\hat\rho(t)=
(e^{\mbox{\boldmath$i$}ta_1},\dots,e^{\mbox{\boldmath$i$}ta_k},1,\dots,1)\in
T^{n+1}\cdot S^1,$$ then $\hat\xi$ is spacelike on
$S^{2n+1,1}-S^1\cdot S^{2(n-k)+1}$.
\end{llemma}

\begin{proof}
Choose an arbitrary point
$v=(z_1,\dots,z_{k},z_{k+1},\dots,z_{n+1},w_1)\in V_0$ such that
\begin{equation}\label{number3}
|z_1|^2+\cdots+|z_k|^2+|z_{k+1}|^2+\cdots+|z_{n+1}|^2-|w_1|^2=0.
\end{equation}

Then
\begin{equation}
\begin{split}
&\rho(t)(z_1,\dots,z_{k},z_{k+1},\dots,z_{n+1},w_1)\\
&=(e^{\mbox{\boldmath$i$}ta_1}z_1,\dots,e^{\mbox{\boldmath$i$}ta_k}z_k,
z_{k+1},\dots,z_{n+1},w_1).
\end{split}\end{equation}
It follows that
$\xi_v=(\mbox{\boldmath$i$}a_1z_1,\dots,\mbox{\boldmath$i$}a_kz_k,0,\dots,0)$
for which
\[\langle\xi_v,\xi_v\rangle=a_1^2|z_1|^2+\cdots+a_k^2|z_k|^2\geq 0.\]
$\langle\xi_v,\xi_v\rangle=0$ if and  only if $z_1=\cdots=z_k=0$. In
this case such a point $v$ satisfies that
\begin{equation}
\begin{split}
P_\RR(v)&=P_\RR((0,\dots,0,z_{k+1},\dots,z_{n+1},w_1))\\
&=\frac
{w_1}{|w_1|}P_\RR((0,\dots,0,\frac{z_{k+1}}{w_1},\dots,\frac{z_{n+1}}{w_1},1))\in
S^1\cdot S^{2(n-k)+1}.
\end{split}\end{equation}
This shows the lemma. (Compare \eqref{compactimage1}.)
\end{proof}

\begin{llemma}\label{lem:D-com}
If the one-parameter group has the form $(${\rm \cf}\ {\bf Case
D}$)$
$$\hat\rho(t)=
(e^{\mbox{\boldmath$i$}ta_1},\dots,e^{\mbox{\boldmath$i$}ta_k},1,\dots,1)e^{\mbox{\boldmath$i$}t}
\in T^{n+1}\cdot S^1,$$ then either one of the following holds.
\begin{enumerate}
\item When all $a_i>0$ or all $a_i<-2$ $(i=1,\dots,k)$, $\hat\xi$ is spacelike on
$S^{2n+1,1}-S^1\cdot S^{2(n-k)+1}$.
\item When all $-2\leq a_i<0$, $\hat\xi$ is timelike on
$S^{2n+1,1}-S^1\cdot S^{2(n-k+\ell)+1}$ for some $\ell< k$.
\item Suppose that there exist $a_i, a_j$ $(1\leq i,j\leq k)$ with
$(i)$ $a_i>0$ or $a_i<-2$, or with $(ii)$ $-2<a_j<0$. Then $\hat
\xi$ is not causal on $S^{2n+1,1}-S^1\cdot S^{2(n-k)+1}$.
\end{enumerate}
\end{llemma}

\begin{proof}
For a point $v$ of \eqref{number3}, it follows similarly as
above:

\[
\xi_v=(\mbox{\boldmath$i$}(a_1+1)z_1,\dots,\mbox{\boldmath$i$}(a_k+1)z_k,
\mbox{\boldmath$i$}z_{k+1},\dots,\mbox{\boldmath$i$}z_{n+1},
\mbox{\boldmath$i$}w_1),
\]
and so
\begin{equation}\label{ecom}
\begin{split}
\langle\xi_v,\xi_v\rangle
 =((a_1+1)^2-1)|z_1|^2+\cdots+((a_k+1)^2-1)|z_k|^2.
\end{split}
\end{equation}

Then the following possibilities occur:\\

\noindent{\bf (1)}\ If all $a_i>0$ or all $a_i<-2$, then
$\langle\xi_v,\xi_v\rangle>0$ and $\langle\xi_v,\xi_v\rangle=0$ if
and only if $z_1=\cdots=z_k=0$ for which $P_\CC(v)\in S^{2(n-k)+1}$.
\vskip0.2cm

\vskip0.2cm \noindent{\bf (2)}\ If all $-2\leq a_i<0$, then
$\langle\xi_v,\xi_v\rangle\leq 0$. Suppose that
$a_{i_1}=\cdots=a_{i_\ell}=-2$ for some $\ell$. By the assumption
that the g.c.m of all $a_i$ is $1$ (\cf Section \ref{sub:compactG}),
note that $\ell<k$. Let
$v=(z_1,\dots,z_{k},z_{k+1},\dots,z_{n+1},w_1)\in V_0$ such that
$z_{i_{\ell+1}}=\cdots=z_{i_{k}}=0$. Then
$\langle\xi_v,\xi_v\rangle=0$ if and only if $P_\CC(v)\in
S^{2(n-k+\ell)+1}$.

\vskip0.2cm \noindent{\bf (3)}\ If there exist $a_i, a_j$ such that
$a_i>0$ (or $a_i<-2$), or $-2<a_j<0$ (\ie $(a_i+1)^2-1>0$,
$(a_j+1)^2-1<0$), then $\langle \xi_v,\xi_v\rangle$ can be taken to
be zero, positive or negative.
\end{proof}
As a consequence,
\par\ \par

\par\noindent{\bf Causality (2)}.\, If $\rho(t)= e^{-\mathrm{i}t}\cdot C_t$ or
 $\rho(t)= e^{\mathrm{i}t}\cdot C_t$  for (i), (ii) of
\eqref{nil-lift}, $\xi$ cannot be lightlike.
\par\ \par

\subsection{Curvature equivalence}\label{curvaturerelation}
 Suppose that $(S^1,M,[g])$ is a
$(2n+2)$-dimensional conformally flat Fefferman-Lorentz parabolic
manifold for which $S^1\leq\mathop{\Conf}_{\rm FLP}(M)$ and let

\begin{equation}\label{devpairS}
(\tilde\rho,\dev):(\tilde S^1,\tilde M)\ra ({\rm
U}(n+1,1)^{\sim},\tilde S^{2n+1,1})
\end{equation}be the developing pair.
Here $\tilde S^1$ is the lift
of $S^1$ to the universal covering $\tilde M$. It is either $S^1$ or
$\RR$. If $\tilde S^1$ induces the vector field $\tilde \xi$
on $\tilde M$, then we note from Definition \ref{def:causalfields}
that
\begin{equation*} \tilde g_x(\tilde \xi,\tilde \xi)=0,\ \ \tilde\xi_x\neq
0 \ (x\in \tilde M).
\end{equation*}
As $\mathop{\dev}(t\cdot x)=\tilde\rho(t)\dev(x)$ ($t\in\tilde
S^1$), $\tilde\rho(\tilde S^1)$ induces the vector field
$\mathop{\dev}_*\tilde \xi$ on the domain ${\dev}(\tilde M)\subset
\tilde S^{2n+1,1}$.
Let $\tilde g$ (respectively $\tilde g^0$) be the lift of $g$
(respectively  the lift of canonical metric $g^0$ to $\tilde
S^{2n+1,1}$). Since $\mathop{\dev}$ is a conformal
immersion, there is a function $u>0$ on $\tilde M$ such that
$\displaystyle u(x)\cdot\tilde g_x(v,w)=\tilde g^0({\dev}_* v,{\dev}_* w)$
In particular,
\begin{equation}\label{conformalmetric}
\tilde g^0_{{\dev}(x)}({\dev}_*\tilde \xi,{\dev}_*\tilde \xi)=0\
  \  ({\dev}(x)\in{\dev}(\tilde M)).
 \end{equation}

\begin{ttheorem}\label{lightlike-th}
Let $(M,[g])$ be a $(2n+2)$-dimensional conformally flat
Fefferman-Lorentz parabolic manifold which admits
$S^1\leq\mathop{\Conf}_{\rm FLP}(M)$. If the $S^1$-action is
lightlike and has no fixed points on $M$, then one of the following
holds.
\begin{itemize}
\item[(i)] $M$ is a Seifert fiber space over a spherical $CR$ orbifold $M/S^1$.
\item[(ii)] The developing pair $(\tilde \rho,\mathop{\dev})$ reduces to
\[
(\mathcal C(\tilde S^1),\tilde M)\to
({\bf R}\times (\mathcal N\rtimes {\rm  U}(n)),\RR\times \mathcal N)
\]where
$\RR\times S^{2n+1}-\RR^1\cdot{\infty}=\RR\times \mathcal N$.
\end{itemize}
\end{ttheorem}

\begin{proof}
Suppose that $\mathsf{G}\neq \{1\}$ for $H=S^1$.
 By the hypothesis, 
there are four possibilities {\bf Cases A, B,C, D}
by Proposition \ref{G-noncompact}. Among them, as
$S^1$ is lightlike, {\bf Causality (1)} and {\bf Causality (2)}
of Section \ref{conf-metric} imply
{\bf Case A} (1), which shows (ii).

When $\mathsf{G}=\{1\}$, first note that $\tilde\rho(\tilde
S^1)={\bf R}$, the center of ${\rm U}(n+1,1)^{\sim}$.
As $\tilde
S^1=\RR$, we have a central extension : $\displaystyle 1\ra \ZZ\ra
\RR\lra S^1\ra 1$. Let $\mathcal Z(\pi)$ be the center of
$\pi=\pi_1(M)$. Since $\tilde S^1$ belongs to the centralizer
$\mathcal Z_{{\rm Diff}(\tilde M)}(\pi)$, it follows $\tilde
S^1\cap\pi\subset\mathcal Z(\pi)$. This shows that $\ZZ=\tilde
S^1\cap\mathcal Z(\pi)=\tilde S^1\cap\pi$. Then this induces a
central extensions:
\begin{equation}\label{group ext}
\begin{CD}
1@>>>  \ZZ @>>> \pi @>>> Q @>>> 1\\
@. \bigcap@. \bigcap@.  {||}@.   \\
1@ >>> \RR @>>> \pi\cdot\RR @>>> Q @>>>1.
\end{CD}
\end{equation}
As $\RR=\tilde S^1$ acts properly and freely on $\tilde M$,
 put $W=\tilde M/\tilde S^1$. Moreover, noting that
$\RR/\ZZ=S^1$ and $\pi$ acts properly discontinuously, the group
$\pi\cdot\RR$ acts properly on $\tilde M$. As a consequence, $Q$
acts properly discontinuously on $W$ with the equivariant fibration:
\begin{equation}\label{Seifert-c}
(\RR,\RR)\ra (\pi\cdot\RR, \tilde M)\lra (Q,W).
\end{equation}
On the other hand, there is the commutative diagram of the holonomy:
\begin{equation}\label{group daigram}
\begin{CD}
{\bf R} @>>> {\rm U}(n+1,1)^{\sim} @>\tilde P>> {\rm PU}(n+1,1)\\
{||}@. \bigcup @. \bigcup@. \\
\tilde\rho(\RR) @>>> \tilde\rho(\pi\cdot\RR) @>\tilde P>>
\hat\rho(Q).
\end{CD}
\end{equation}
Then the developing pair \eqref{devpairS} induces an
equivariant developing map
on the quotient space:
\[
(\hat\rho,\hat{\dev})\colon (Q, W)\ra ({\rm PU}(n+1,1), S^{2n+1}),
\]where $S^{2n+1}=\tilde S^{2n+1,1}/{\bf R}$.
Since $({\rm PU}(n+1,1), S^{2n+1})$ is the spherical $CR$-geometry,
$W$ inherits a spherical $CR$-structure on which $Q$ acts as
$CR$-transformations. Taking the quotient of \eqref{Seifert-c},
$S^1\ra M\ra M/S^1$ is a Seifert fiber space over the $CR$-orbifold
$M/S^1=Q\backslash W$.
\end{proof}

Using Theorem \ref{lightlike-th} we can prove the following
equivalence.
\begin{ttheorem}\label{equivalence-flatness}
Let $(S^1\times N,g)$ be a Fefferman-Lorentz manifold for a strictly
pseudoconvex $CR$-manifold $(N,(\om_N,J_N))$ of dimension $2n+1\geq
3$ in which $g=\sigma\odot P^*\om_N+d\om_N(J_NP_*-,P_*-)$ is a
Fefferman metric. Then $(S^1\times N,g)$ is conformally flat if and
only if $(N,(\om_N,J_N))$ is spherical $CR$.
\end{ttheorem}

\begin{proof}
 By Theorem \ref{lightlike-th}, the case (i) or (ii) occurs. If (i) occurs, then
$(Q,W)=(\pi_1(N),\tilde N)$ with $\tilde S^1=\RR$ for which
there is
a developing map
\begin{equation}\label{deve-cr}
(\hat\rho,\hat{\dev}):(Q, \tilde N)\ra ({\rm PU}(n+1,1), S^{2n+1}).
\end{equation}
We have to check that the  spherical $CR$-structure $(\om,\hat J)$
induced by $\hat{\dev}$ coincides with the original one
$(\om_N,J_N)$ on $N$. The contact form $\om$ is obtained as
$$P^*\om=g(\mathcal S, -)$$ from
\eqref{contact} of Proposition \ref{1G-noncompact} and the complex
structure $\hat J$ is defined by
 $$\hat{\dev}_*\hat J=J_0\hat{\dev}_*$$ on $\mathop{\Ker}\, \om$
from \eqref{complexCR}.

Let $\mathcal S$ be the vector field induced by $S^1$ on $S^1\times
N$. Since $\displaystyle \sigma(\mathcal S)=\frac 1{n+2}$ from
\eqref{sigvavalue}, it follows that $\displaystyle g(\mathcal
S,-)=\frac 1{n+2}P^*\om_N(-)$ and so

\begin{equation}\label{equlitysig}
 \om_N=(n+2)\om.
 \end{equation}

If we note that $S^1$ acts as lightlike isometries of $(S^1\times
N,g)$, then it satisfies also Proposition \ref{1G-noncompact}
(\cf Remark \ref{g-invariant}).
Then from \eqref{crequiv},

$$(n+2)\hat u\cdot \om=\hat{\dev^*}\om_0,$$
which implies that $\hat{\dev_*}(\mathop{\Ker \, \om})=
\mathop{\Ker\, \om_0}$.\\

 Let $\mathop{\dev}^*g^0=u\cdot g$ as
before. If $P_*X,P_*Y\in \mathop{\Ker \, \om}$, then
\[
g(X,Y)=\sigma\odot
P^*\om_N(X,Y)+d\om_N(J_NP_*X,P_*Y)=d\om_N(J_NP_*X,P_*Y)\]  Noting
$P^*\hat u=u$,
\begin{equation}\label{JN}
u\cdot g(X,Y)=\hat u\cdot d\om_N(J_NP_*X,P_*Y).
\end{equation}

On the other hand, using
$\hat{\dev}_*\hat J=J_0\hat{\dev}_*$ with \eqref{downCR},
\begin{equation*}
\begin{split}
{\dev}^*g^0(X,Y)&=g^0({\dev}_*X,{\dev}_*Y)\\
&=d\om_0(J_0P_*({\dev}_*X),P_*({\dev}_*Y)\\
&=d\om_0(J_0\hat{\dev}_*P_*X,\hat{\dev}_*P_*Y)\\
&=d\om_0(\hat{\dev}_*\hat JP_*X,\hat{\dev}_*P_*Y)\\
&=\hat{\dev^*}d\om_0(\hat JP_*X,P_*Y).
\end{split}\end{equation*}
 Noting that $(n+2)\hat u\cdot d\om=\hat{\dev^*}d\om_0$ on
$\mathop{\Ker \, \om}$, it follows by \eqref{equlitysig} that
\begin{equation}\label{OriJ}\begin{split}
{\dev}^*g^0(X,Y)&=(n+2)\hat u\cdot d\om(\hat JP_*X,P_*Y)\\
&=\hat u\cdot  d\om_N(\hat JP_*X,P_*Y).
\end{split}\end{equation} Compared \eqref{JN} and  \eqref{OriJ}, we
conclude that \[\hat J= J_N.\] Hence $(\mathop{\Ker}\, \om,\hat J)=(
\mathop{\Ker}\, \om_N, J_N)$ so that $(N,(\om_N,J_N))$ is a
spherical $CR$-manifold.
\\

We have to show that the case (ii) of Theorem \ref{lightlike-th}
does not occur.  If (ii) occurs, then we have a developing pair
by Proposition \ref{G-noncompact}:
\begin{equation}\label{devpairN} (\tilde\rho,\dev):(\tilde
S^1,\RR\times \tilde N)\ra
({\bf R}\times (\mathcal N\rtimes {\rm  U}(n)), \RR\times \mathcal N).
\end{equation}
Here $\tilde S^1=\RR$.
Let $\tilde{\mathcal S}$ be the vector field induced by $\RR$ on $\RR\times \tilde N$
as before. Put $\mathop{\dev}_*\tilde{\mathcal S}=\tilde{\mathcal S}'$.
Let
$$g_\mathcal N=\sigma_\mathcal N\odot P^*\om_\mathcal
N+d\om_\mathcal N(JP_*-,P_*-)$$ be the Lorentz metric on $\RR\times
\mathcal N$ which is conformal to the standard metric $g^0$.
(Compare Section \ref{conf-metric}.) Note that
 $g_\mathcal N(\tilde{\mathcal S}',\tilde{\mathcal S}')=0$ because
$\tilde{\mathcal S}$ is lightlike. In this case, {\bf Causality (1)}
shows that $\tilde{\mathcal S}'$ is
the characteristic vector field, \ie $\om_\mathcal N(P_*\tilde{\mathcal S}')=1$.
Moreover, Proposition \ref{G-noncompact}
with {\bf Causality (2)} implies that 
the lightlike vector field $\tilde{\mathcal S}'$ is
 of type ${\rm (i)}$ of \eqref{nil-lift}.
  As $\displaystyle \sigma_\mathcal N=\frac 1{n+2}dt$
(\cf Section \ref{nipotentcase}), it follows that
 $\sigma_\mathcal N(\tilde{\mathcal S}')=0$.

As before, there exists a function $u>0$ on $\RR\times\mathcal N$
such that $\mathop{\dev}^*g_\mathcal N=u\cdot g$.
 Noting that
$P_*\tilde{\mathcal S}'$ is characteristic, a calculation shows that
\begin{equation*}
\begin{split}
g_\mathcal N({\dev}_*\tilde{\mathcal S},{\dev}_*V)&=
g_\mathcal N(\tilde{\mathcal S}',{\dev}_*V)\\
&=\sigma_\mathcal N\odot P^*\om_\mathcal N(\tilde{\mathcal S}',
{\dev}_*V)
+ d\om_\mathcal N(JP_*\tilde{\mathcal S}',P_*{\dev}_*V)\\
&=\sigma_\mathcal N({\dev}_*V) \\
 u\cdot g(\tilde{\mathcal S},V)&=
u(\sigma\odot P^*\om_N(\tilde{\mathcal S},V)
+d\om_N(JP_*\tilde{\mathcal S},P_*V))\\
&=\frac u{n+2}\cdot\om_N(P_*V)\, \ (P_*\tilde{\mathcal S}=0).
\end{split}
\end{equation*}It follows that
\begin{equation}
{\dev}^*\sigma_\mathcal N=\frac u{n+2}\cdot P^*\om_N.
\end{equation}
It is easy to see that $\displaystyle \frac
{u^{n+1}}{(n+2)^{n+1}}P^*(\om_N\we (d\om_N)^{n})=
{\dev}^*(\sigma_\mathcal N\we (d\sigma_\mathcal N)^n)$ on $\RR\times
\tilde N$. As $d\sigma_\mathcal N=0$ as above, it follows that
$P^*(\om_N\we(d\om_N)^{n})=0$ so that $\om_N\we (d\om_N)^{n}=0$ on
$\tilde N$, which contradicts that $\om_N$ is a contact form on
$\tilde N$. Therefore the case (ii) of
Theorem \ref{lightlike-th} cannot occur. This proves the necessary condition.\\

Suppose that $N$ is spherical $CR$. There exists a collection of
charts
 $\{U_\al,\varphi_\al\}_{\al\in\Lambda}$
such that $\varphi_\al :U_\al\ra \varphi_\al(U_\al)\subset S^{2n+1}$
is a homeomorphism. Consider the pullback of the $S^1$-bundle:
\begin{equation}\label{eq:pullback}
\begin{CD}
S^1@>>> S^1\\
@VVV @VVV \\
S^1\times U_\al@>\tilde \varphi_\al>> S^1\times S^{2n+1}\\
@VVV @V{P_\CC}VV \\
U_\al@>\varphi_\al>> S^{2n+1}
\end{CD}
\end{equation}in which
\begin{equation}\label{exten-h}
\tilde \varphi_\al(t,x)=(t,\varphi_\al(x))
\end{equation}
When $U_\al\cap U_\be\neq \emptyset$, the local change
$\varphi_\al\circ \varphi_\be^{-1}$ extends to an automorphism $h\in
{\rm PU}(n+1,1)$ of $S^{2n+1}$. Put $U=\varphi_\be(U_\al\cap
U_\be)\subset S^{2n+1}$. Consider the local diffeomorphism:
\begin{equation}\label{localglobal}
\tilde h=\tilde\varphi_\al\circ \tilde\varphi_\be^{-1}\ : S^1\times
U \ra S^1\times S^{2n+1}
\end{equation}for which
\begin{equation}\label{explicith}
\tilde h(t,z)=(t,hz).
\end{equation}
If we note that ${\mathrm U}(n+1,1)$ acts invariantly on $V_0$
such that $S^1\times S^{2n+1}\subset V_0$,
then there exists an element $f\in \mathrm{U}(n+1,1)$ with $Pf=h$
which satisfies that
\begin{equation}\label{localconf}
\tilde h=f| S^1\times U.
\end{equation}
By Proposition \ref{confflatlo},
${\mathrm{U}}(n+1,1)$ acts conformally on $S^1\times
S^{2n+1}$ with respect to $g^0$ so
it follows that
\[
\tilde h^*g^0={f}^*g^0=v\cdot g^0 \ \, (\exists\,v>0).\] Since
$\tilde h$ is a local conformal diffeomorphism, $\tilde h$ extends
to a global conformal transformation of $S^1\times S^{2n+1}$ by the
Liouville's theorem. By uniqueness,
\[
\tilde h=f\ \, \mbox{on}\ S^1\times S^{2n+1}.
\]
As a consequence, the local change $\tilde \varphi_\al\circ \tilde
\varphi_\be^{-1}$ extends to an automorphism $\tilde h\in {\rm
U}(n+1,1)$ of $S^{2n+1}\times S^1$. Therefore the charts
$\{U_\al\times S^1 ,\tilde \varphi_\al\}_{\al\in\Lambda}$ of
$N\times S^1$ gives a uniformization with respect to $({\rm
U}(n+1,1), S^{2n+1}\times S^1)$. As $({\rm U}(n+1,1), S^{2n+1}\times
S^1)$ is the lift of $(\hat {\rm U}(n+1,1), S^{2n+1,1})$, $N\times
S^1$ is a conformally flat Fefferman-Lorentz manifold.

\end{proof}

\section{Application to Obata \& Ferrand's theorem}
\label{subsec:Obata-Ferrand}

\subsection{Noncompact conformal group actions} Let $\mathcal C$ be a
closed noncompact subgroup of $\mathop{\Diff}(M)$. Suppose that
$\mathcal C$ acts \emph{analytically} on $M$. Contrary to compact
group actions, noncompact (analytic) Lie group actions on a compact
manifold is quite different. For example, there is a noncompact
analytic action $(\mathcal C,M)$ such that the set of nonprincipal
orbits $M_0=\{x\in M\, |\, {\rm dim}\, \mathcal C\cdot x<{\rm dim}\,
\mathcal C\}$ coincides with $M$. (See \cite{ri},  \cite[3.1
Theorem, p.79]{bre} for instance.) We give a sufficient condition
 that $M_0$ is nowhere dense in $M$
for noncompact Lorentz groups $\mathcal C$ acting analytically
on a Lorentz manifold $(M,\mathsf{g})$. Let $\mathcal C=S^1\times
\RR$ which is closed in $\mathop{\Diff}(M)$. If $S^1$ does not act
as Lorentz isometries with respect to $\mathsf{g}$, then we put
\begin{equation}\label{Smetric}
\tilde{{\mathsf g}}=\int_{S^1}^{} h^*{\mathsf g} dh
\end{equation} where $ds$ is a right-invariant Haar measure on
$S^1$. Then $S^1$ acts as Lorentz isometries with respect to $\tilde
{\mathsf g}$. On the other hand, as $S^1$ acts conformally with
respect ${\mathsf g}$ , $h^*{\mathsf g}=\lam_h\cdot {\mathsf g}$ for
some function $\lam_h>0$ on $M$. Letting $\displaystyle
\tau(x)=\mathop{\int}_{S^1}^{}\lam_h(x)dh$ $(x\in M)$, it follows
that $\tilde{\mathsf g}=\tau\cdot{\mathsf g}$ on $M$. We obtain a
Lorentz metric $\tilde {\mathsf g}$ conformal to ${\mathsf g}$. So
if $\mathcal C=S^1\times \RR$ acts conformally on $(M,g)$, we may
assume that $S^1$ acts as isometries within the conformal class of
the Lorentz metric $g$.

\begin{ppro}\label{sameorbit}
If a Lorentz $(n+2)$-manifold $(M,g)$ admits a closed two
dimensional subgroup $\mathcal C$ isomorphic to $S^1\times \RR$
where $S^1$ consists of lightlike conformal transformations, then
\[ V=\{x\in M\,|\, {\rm dim}\, \mathcal C\cdot x=2\} \]is a dense
open subset of $M$.
\end{ppro}

\begin{proof}
We suppose that $S^1$ acts isometries.
 Let $\tilde F$ (respectively $F$) be the fixed point set
of $\RR$ (respectively $S^1$). Note that $S^1$ leaves ${\tilde F}$
invariant.
If $E$ is the set of exceptional orbits of $S^1$, then
$S^1$ acts freely on the complement $M_0=M-{(E\cup F\cup \tilde
F)}$. Note that $M_0$ is a dense open subset of $M$. There is a
principal bundle over the orbit space $N_0=M_0/S^1$;
\begin{equation}\label{princi-action}
\begin{CD}
S^1@>>> M_0@>P>>N_0.
\end{CD}
\end{equation}
Suppose that ${\rm dim}\, \mathcal C\cdot x=1$ for some open subset
$U$ of $M_0$ $(\forall x\, \in U)$. If we set $\{\varphi_t\}_{t\in
\RR}=\RR$, then it follows that
\[
\mathcal C\cdot x=\{\varphi_tx\}_{t\in \RR}=S^1\cdot x.\] Let
$\xi_x$ be the vector induced by $\mathcal C\cdot x$ $(x\in U)$. By
the hypothesis, $\xi$ is a lightlike (Killing) vector field on $U$.
For an arbitrary point $y\in U$, as $\varphi_t y\in S^1\cdot y$,
there exists an element $h^y_t\in S^1$ such that
\begin{equation}\label{equaphi}
\varphi_t y= h^y_t\cdot y.
\end{equation}
This implies that $P\circ \varphi_t=P$ on $M_0$. Put $z=\varphi_t y
=h_t y$ where we let $h_t=h_t^y$ for brevity.
For a vector $v_y\in
T_yM_0$, we have that
 $P_*\varphi_{t*}v_y=P_*v_y=P_*h_{t*}v_y$. Since
$\varphi_{t*}v_y,\, h_{t*}v_y\in T_zM_0$, it follows that
\[\varphi_{t*}v_y=h_{t*}v_y+a\xi_z \ \,(\exists\, a\in \RR).\]
As $\xi_y$ is lightlike, we can find a vector $\eta_y\in T_yM_0$
such that
\begin{equation*}
g(\eta_y,\eta_y)=0,\ \, g(\xi_y,\eta_y)=1.
\end{equation*}
As above, there exists an element $b\in \RR$ such that
\[ \varphi_{t*}\eta_y=h_{t*}\eta_y+b\xi_z.\]
Since $g(\varphi_{t*}\eta_y,\varphi_{t*}\eta_y)=\lam_t(y)\cdot
g(\eta_y,\eta_y)=0$ and $h_t y=z$, a calculation shows that
\begin{equation}\begin{split}
0&=g(h_{t*}\eta_y+b\xi_z,h_{t*}\eta_y+b\xi_z)\\
&=2bg(h_{t*}\eta_y,\xi_z)=2bg(h_{t*}\eta_y,h_{t*}\xi_y)\\
&=2bg(\eta_y,\xi_y)=2b,\\
\end{split}
\end{equation}so it follows that $\varphi_{t*}\eta_y=h_{t*}\eta_y$.

Noting that $\{ \xi_y,\eta_y\}$ spans a nondegenerate plane of
signature $(1,1)$, there exists a vector $v_y$ such that
$g(v_y,v_y)=1,\, g(\xi_y,v_y)=g(\eta_y,v_y)=0$. There are
$n$-independent such vectors. The set of those vectors
 with $\{ \xi_y,\eta_y\}$ constitutes
$T_yM_0$. As above, let $\varphi_{t*}v_y=h_{t*}v_y+a\xi_z$.
Similarly, using $\varphi_{t*}\eta_y=h_{t*}\eta_y$, the equation
$g(\varphi_{t*}\eta_y,\varphi_{t*}v_y)=0$ shows that $a=0$, \ie
$\varphi_{t*}v_y=h_{t*}v_y$. From these calculations, we obtain that
\begin{equation}\label{all*}
\varphi_{t*}X_y=h_{t*}X_y,\  \, (\forall\, X_y\in T_yM_0).
\end{equation}
Now, noting $h_t\in S^1$,

\begin{equation}\label{isom}
\varphi_t^*g(X_y,Y_y)=g(h_{t*}X_y,h_{t*}Y_y)=g( X_y,Y_y)\ \,
(\forall X_y,Y_y\in T_yM_0).
\end{equation}

On the other hand, since $\RR$ acts conformally, there exists a
positive function $\lam_t$ on $M$ such that
$\varphi_t^*g=\lam_t\cdot g$ for each $t$, \eqref{isom} implies that
$\lam_t(y)=1$. This is true for an arbitrary point $y\in U$, so
$\lam_t=1$ on $U$. In particular, $\varphi_t$ $(\forall\,t\in \RR)$
becomes a Lorentz isometry on $U$ (and so is on $M$ by analyticity).

Recall from \eqref{equaphi} that $\varphi_t y=h_t\cdot y$. Since
$\RR= \{\varphi_t\}_{t\in \RR}$, there exists an element $a\in \RR$
such that $\varphi_a y=y$. (In fact, if $\displaystyle
a(y)=\mathop{\min}_{t\in \RR^+}^{}\{t\, |\ \varphi_t y=y\}$, we put
$a=a(y)$.) Then it follows that $y=\varphi_a y=h_a y$.  As $S^1$
acts freely on $M_0$, $h_a=1$. From \eqref{all*}, we have that
\begin{equation}\label{all**}
\varphi_{a*}X_y=X_y\  \, (\forall\, X_y\in T_yM_0).
\end{equation}
Since $\varphi_a$ is a Lorentz isometry, if $\ga$ is any geodesic
issuing from $y$, then $\varphi_a\ga$ is also a geodesic on $U$.
From \eqref{all**}, the uniqueness of geodesic implies that
$\varphi_a\ga=\ga$ on $U$. Hence $\varphi_a={\rm id}$ on $U$. By
analyticity, $\varphi_a={\rm id}$ on $M$.  Letting $\ZZ=\langle
na\rangle_{n\in \ZZ}$ so that $S^1=\RR/\ZZ$, $\mathcal C$ would be
isomorphic to $S^1\times S^1$. This contradicts our hypothesis that
$\mathcal C$ is noncompact. Hence the subset $\{x\in M\ |\ {\rm
dim}\, \mathcal C\cdot x=2\}$ is dense open in $M$.

\end{proof}

\begin{ttheorem}\label{th:o-l} Let $M=S^1\times N$ be a compact Fefferman-Lorentz
manifold and $\mathcal C_{{\Conf}(M,g)}(S^1)$ the centralizer of
$S^1$ in ${\Conf}(M,g)$. Suppose that
  $\mathcal C_{{\Conf}(M,g)}(S^1)$ contains a closed
noncompact subgroup of dimension $1$ at least. Then $M$ is
conformally equivalent to the two-fold cover $S^1\times S^{2n+1}$ of
the standard Lorentz manifold $S^{2n+1,1}$.
\end{ttheorem}

\begin{proof}
We can choose a closed subgroup $\mathcal C=S^1\times \RR$ from
$\mathcal C_{{\Conf}(M,g)}(S^1)$ by the hypothesis. Here recall that
the vector field $\mathcal S$ generated by $S^1$ of $M=S^1\times N$
is lightlike.

Recall from \eqref{f-lorentz} that
\begin{equation}\label{lorentmetricf}
g=\sigma\odot P^*\om+d\om(JP_*-,P_*-)
\end{equation}is a Lorentz metric on a Fefferman-Lorentz manifold $M=S^1\times N$ where
$P\colon S^1\times N\to N$ is the projection. Then $\mathcal C$
induces an action of $\RR$ on the quotient $N$ such that $P$ is
equivariant:
\begin{equation}\label{projection}
 P:(\mathcal C,M)\ra (\RR,N).
 \end{equation}
If $\{\varphi_t\}_{t\in \RR}$ is a $1$-parameter group of $\RR$ of
$\mathcal C$, then there exists a $1$-parameter group of
$\{\hat\varphi_t\}_{t\in \RR}$ such that
\[
P\circ\varphi_t=\hat\varphi_t\circ P.\]

Since $\RR$ acts as conformal transformations with respect to $g$,
there exists a function $\lam_t:M\ra \RR^+$ such that
\begin{equation}\label{conformalmap}
\varphi_t^*g=\lam_t\cdot g.
\end{equation}If $h\in S^1$, since
$h^*\varphi_t^*g=\varphi_t^*h^*g=\varphi_t^*g$ and
$h^*\varphi_t^*g=h^*(\lam_t\cdot g)=h^*\lam_t\cdot g$, it follows
that $h^*\lam_t=\lam_t$  $(\forall\,h\in S^1)$. So $\lam_t$ factors
through a function $\hat\lam_t: N\ra \RR^+$ $(\forall\, t\in \RR)$.
We note also that $\varphi_{t*}\mathcal S=\mathcal S$ and
$P_*\mathcal S=0$. Then
\begin{equation}\label{equality1}\begin{split}
\varphi_t^*g(X,\mathcal S)&=\frac 1{n+2}\om(P_*\varphi_*X)=\frac
1{n+2}P^*\hat\varphi*\om(X)\\
&=\lam_t\cdot g(X,\mathcal S)=\frac 1{n+2}P^*\hat\lam_t\cdot
P^*\om(X). \end{split}\end{equation} it follows that
\begin{equation}\label{conformlaomega}
\hat\varphi_t^*\om=\hat\lam_t\cdot \om\ \ (\forall\, t\in \RR).
\end{equation}
This implies that
\begin{equation}\label{prseevKer}
\hat\varphi_{t*}{\rm Ker}\,\om={\rm Ker}\, \om. \end{equation}
Recall that there is a complex structure $J$ on ${\rm Ker}\, \om$.
For convenience, we put $\tilde J$ on $P^*{\rm Ker}\, \om$ formally
such that $P_*:P^*{\rm Ker}\, \om\ra {\rm Ker}\, \om$ is almost
complex, \ie $P_*\circ \tilde J=J\circ P_*$.

Let $X,Y\in P^*{\rm Ker}\, \om$. Calculate
\begin{equation}\label{g-iso1}
\begin{split}
\varphi_{t}^{*}g(-X,Y)&=g(-\varphi_{t*}X,\varphi_{t*}Y)\\
&=d\om(JP_*(-\varphi_{t*}X),P_*(\varphi_{t*}Y)) \ \
\text{by}\ \eqref{lorentmetricf}\\
&=d\om(-J\hat\varphi_{t*}P_*X,\hat\varphi_{t*}P_*Y)\\
&=d\om(\hat\varphi_{t*}P_*X,J\hat\varphi_{t*}P_*Y).\\
\end{split}
\end{equation}

 Noting that $\hat\lam_t\cdot d\om=\hat\varphi_{t}^*d\om $ on ${\rm Ker}\, \om$,
 calculate
\begin{equation}\label{g-iso2}
\begin{split}
&\varphi_{t}^{*}g(-X,Y)=\lam_t\cdot g(-X,Y)\\
&\ \ =\lam_t\cdot d\om(JP_*(-X),P_*Y)=\lam_t\cdot d\om(P_*X,JP_*Y)\\
&\ \ =\lam_t\cdot d\om(P_*X,P_*(\tilde J)Y)=P^*(\hat\lam_t\cdot d\om)(X,\tilde JY)\\
&\ \ =\hat\lam_t\cdot d\om(P_*X,P_*(\tilde JY))=\hat\lam_t\cdot d\om(P_*X,JP_*Y)\\
&\ \ =\hat\varphi_{t}^*d\om(P_*X,JP_*Y)\\
&\ \ =d\om(\hat\varphi_{t}^*P_*X,\hat\varphi_{t}^*JP_*Y).
\end{split}
\end{equation}As $d\om$ is nondegenerate on ${\rm Ker}\, \om$, we
conclude that
\begin{equation}\label{holoJ}
\hat\varphi_{t*}J=J\hat\varphi_{t*}\ \ \text{on}\ {\rm Ker}\, \om.
\end{equation}Let
${\rm Aut}_{CR}(N)$ be the group of $CR$-transformations of
$(\om,J)$ on $N$. By the definition, $\{\hat\varphi_{t}\}_{t\in
\RR}\subset {\rm Aut}_{CR}(N)$. Note that $\{\hat\varphi_{t}\}_{t\in
\RR}$ is closed by the hypothesis. Moreover, the action of
$\{\hat\varphi_{t}\}_{t\in \RR}$ is nontrivial on $N$ by Proposition
\ref{sameorbit}. Hence, $\{\hat\varphi_{t}\}_{t\in \RR}$ is a closed
noncompact subgroup in ${\rm Aut}_{CR}(N)$. It follows from the
$CR$-analogue of Obata-Ferrand rigidity (for example \cite{kam3},
\cite{we}, \cite{cf}, \cite{lee1}, \cite{sh-cr}) that $N$ is
$CR$-isomorphic to the standard sphere $S^{2n+1}$.
 Then $(S^1\times N,g)$ is conformally flat by Theorem
 \ref{equivalence-flatness}.
Let $\tilde {\mathcal C}$ be a lift of $\mathcal C$ to $\tilde M$.
By Proposition \ref{pro:flat-almost complex}, we have the developing
pair:
\begin{equation}\label{LorentzSeifert}
(\tilde\rho,\dev): (\tilde {\mathcal C}, \tilde M)\lra ({\rm
U}(n+1,1)^{\sim},\tilde S^{2n+1,1}).
\end{equation}
Recall that ${\bf R}$ is the center of ${\rm U}(n+1,1)^{\sim}$ which
is the kernel of projection $\tilde P\colon {\rm U}(n+1,1)^{\sim}\to
{\rm PU}(n+1,1)$ and $\mathcal R$ is the center of the Heisenberg
group $\mathcal N$ in ${\rm PU}(n+1,1)$ (\cf \eqref{groupgram}). Let
$\tilde S^1$ be a lift of  lightlike one-parameter subgroup to
$\tilde {\mathcal C}$. We show that

\begin{equation}\label{eq:1-dimflow}
\tilde \rho:{\tilde S}^1\lra {\bf R}
\end{equation}is isomorphic.
For this, put $\displaystyle \mathsf{G}=\overline{\tilde P\circ
\tilde\rho(\tilde S^1)}\subset{\rm PU}(n+1,1)$ as in \eqref{cloG}.
{\bf Causality ${\bf (1)}$} and {${\bf (2)}$} in Section
\ref{Fefferman-Lorentz-CR manifolds} yield that
\begin{equation*}
\tilde\rho(\tilde H_1)= \left\{\begin{array}{lr} \mathcal R&\,
\mbox{if}\,
\, \mathsf{G}\neq \{1\},\\
{\bf R}&\, \ \, \mbox{if}\ \, \mathsf{G}=\{1\}.\\
\end{array}\right.
\end{equation*}If $\mathsf{G}\neq\{1\}$, then
by (ii) of Theorem \ref{lightlike-th} the developing pair reduces to
\[
(\tilde \rho,\mathop{\dev}): (\mathcal C(\tilde S^1),\tilde M)\to
({\bf R}\times (\mathcal N\rtimes {\rm  U}(n)),\RR\times \mathcal N)
\]where $\RR\times S^{2n+1}-\RR^1\cdot{\infty}=\RR\times \mathcal N$. Since
$\tilde {\mathcal C}\subset \mathcal C(\tilde S^1)$ and
 ${\bf R}\times(\mathcal N\rtimes {\rm U}(n))$ is transitive on
${\bf R}\times\mathcal N$ with compact stabilizer $\rm {U}(n)$,
$\tilde M$ admits a $\pi\cdot\tilde {\mathcal C}$-invariant
\emph{Riemannian metric} $\tilde{\mathsf{g}}$. Taking the quotient,
it follows that $\mathcal C\leq {\rm Isom}(M,\mathsf{g})$. As $M$ is
compact, ${\rm Isom}(M)$ is compact. Since $\mathcal C$ is a closed
noncompact subgroup (in $\mathop{\Diff}(M)$) by our hypothesis, this
case cannot occur.\\

Then $\mathsf{G}=\{1\}$ and so $\tilde\rho(\tilde S^1)={\bf R}$. It
follows from (i) of Theorem \ref{lightlike-th} that $M$ is a Seifert
fiber space over a spherical $CR$ orbifold $M/S^1$. In our case,
$M/S^1=\tilde M/\tilde S^1=N$ which is simply connected. We obtain
the following commutative diagram:
\begin{equation}\label{verydiagarm}
\begin{CD}
\tilde S^1@>\tilde \rho>> {\bf R}\\
@VVV @VVV \\
\tilde M@>\tilde {\dev}>> {\bf R}\times S^{2n+1}\\
@V{P}VV @V{\tilde P}VV \\
N@>\hat\dev>> S^{2n+1}.
\end{CD}
\end{equation}
Moreover, the $CR$-structure $(\om,J)$ on $N$ coincides with the
pullback of the standard $CR$-structure $(\om_0,J_0)$ of $S^{2n+1}$.
In fact, we have shown in the proof of Theorem
\ref{equivalence-flatness} that
\begin{equation}\begin{split}
 \hat u\cdot \om&=\hat{\dev^*}\om_0,\\
 \hat{\dev}_*\hat J&=J_0{\dev}_* \ \text{ on}\ \mathop{\Ker}\, \om
\end{split}\end{equation}
As $(N, (\om,J))$ is $CR$-isomorphic to $S^{2n+1}$ as above, there
exists an element $h\in {\rm PU}(n+1,1)$ such that
$h\circ\hat{\dev}:N\ra S^{2n+1}$ is a $CR$-diffeomorphism. (As a
consequence, $\hat{\dev}$ itself is a $CR$-diffeomorphism.) By the
diagram \eqref{verydiagarm}, $\tilde {\dev}:\tilde M\ra {\bf
R}\times S^{2n+1}$ is a conformal diffeomorphism. Taking a quotient,
$\tilde {\dev}$ induces a conformal diffeomorphism
 ${\dev}: M\ra S^1\times S^{2n+1}$.

\end{proof}

\begin{rremark}\label{rem:1-theorem8.12}
Let ${\rm SL}(2,\RR)/\Gamma$ be a Lorentz space form of negative
constant curvature where $\Gamma\leq {\rm O}(2,2)^0={\rm
SL}(2,\RR)\cdot{\rm SL}(2,\RR)$ is a subgroup acting properly
discontinuously on ${\rm SL}(2,\RR)$ $($\cf \cite{kul-ray}$)$. If we
choose
$$\Gamma\leq {\rm SO}(2)\mathop{\times}_{\ZZ_2}^{} {\rm
SL}(2,\RR),$$ then ${\rm SL}(2,\RR)/\Gamma$ is a spherical
$CR$-space form because there is a canonical identification:
\begin{equation*}\begin{split}
& {\rm SO}(2)\mathop{\times}_{\ZZ_2}^{} {\rm SL}(2,\RR)={\rm
U}(1,1) < {\rm U}(1)\times {\rm U}(1,1) < {\rm U}(2,1),\\
& {\rm SL}(2,\RR)=S^3-S^1. \end{split}
 \end{equation*}

{\rm {\bf (1)}}\ Thus $M=S^1\times {\rm SL}(2,\RR)/\Gamma$ is a
conformally flat Fefferman-Lorentz $4$-manifold on which $S^1$ acts
as lightlike isometries by the definition.

 {\rm {\bf (2)}}\  Recall that $(\rm {PO}(4,2),S^{3,1})$ is the
conformally flat Lorentz geometry. Let $(\rm {O}(4,2),S^1\times
S^{3})$ be the two fold covering. The subgroup of $\rm {O}(4,2)$
preserving $S^1\times (S^{3}-S^1)=S^1\times {\rm SL}(2,\RR)$ is
isomorphic to $\rm {O}(2)\times \mathrm {O}(2,2)$. When we restrict
$\Gamma$ to $\{1\}\times\mathrm {SL}(2,\RR)$, we obtain a
conformally flat Lorentz parabolic manifold $M=S^1\times \mathrm
{SL}(2,\RR)/\Gamma$ on which $S^1={\rm SO}(2)$ acts as spacelike
isometries. The subgroup $$\mathrm {SO}(2)\times (\mathrm
{SL}(2,\RR)\times \{1\})$$ of $\mathrm {O}(2)\times \mathrm
{O}(2,2)$ acts conformally on $M$. Let $\mathrm {SL}(2,\RR)=KAN$ be
the decomposition as usual. In particular,

{\rm {\bf (3)}}\ $M$ admits the two-dimensional closed noncompact
conformal group $C=S^1\times N\times \{1\}$ consisting of spacelike
and lightlike transformations. However, $C$ does not belong to ${\rm
U}(2,1)$ because
\[\mathrm U(1)=C\cap{\rm U}(2,1)\leq  (\mathrm {O}(2)\times \mathrm {O}(2,2))\cap {\rm U}(2,1)=
\mathrm U(1)\times {\rm U}(1,1)\] so $C$ does not preserve the
Fefferman-Lorentz parabolic structure.
\end{rremark}

\end{document}